\newif\ifxivarxive
\xivarxivetrue 

\ifxivarxive
    \documentclass[onefignum,onetabnum]{siamart171218}
    \newcommand{\sep}{\text{, }}
\else
    \documentclass[review,onefignum,onetabnum]{siamart171218}
    \newcommand{\sep}{\text{, }}
\fi

\usepackage{amsmath,amssymb,amsfonts,graphicx,hyperref}
\usepackage[margin=1in]{geometry}

\usepackage{amstext}
\usepackage{stmaryrd}
\usepackage{algpseudocode}
\usepackage[algo2e,linesnumbered,vlined,commentsnumbered,ruled]{algorithm2e}

\usepackage[normalem]{ulem}
\usepackage{graphicx}

\usepackage{array}
\usepackage{verbatim}
\usepackage{booktabs}
\usepackage{calc}
\usepackage{textcomp}
\usepackage{multirow}

\usepackage[all]{xy}

\usepackage{cancel}
\usepackage{mathtools}
\usepackage{placeins}
\usepackage{xspace}
\usepackage[colorinlistoftodos,prependcaption,textsize=small]{todonotes}

\usepackage{tikz,tikzscale}
\usetikzlibrary{calc}
\usetikzlibrary{shadings}

\usepackage[mode=multiuser,status=draft]{fixme} 
\fxsetup{envlayout=color}
\fxsetup{innerlayout=inline}
\fxsetup{targetlayout=colorcb}
\FXProvidesTargetLayout{color}

\fxusetheme{color}
\FXRegisterAuthor{ss}{envss}{SS}
\FXRegisterAuthor{mr}{envmr}{MR}
\FXRegisterAuthor{mw}{envmw}{MW}
\FXRegisterAuthor{mk}{envmk}{MK}

\usepackage{subcaption}
\usepackage[format=plain,indention=.5cm]{caption}
\usepackage{pgfplots}

\usepackage{bm,upgreek}

\usepackage{graphicx} 
\usepackage[linesnumbered,vlined,commentsnumbered,ruled]{algorithm2e}

\definecolor{sbmgreen}{RGB}{0,100,0}

\newcommand{\plotlinestyleone}[1][solid]{%
    \addplot[color=orange,mark options={solid}, mark=x, mark size=2pt,line
        width=1pt,#1]%
}
\newcommand{\plotlinestyletwo}[1][solid]{%
    \addplot[color=sbmgreen,mark options={solid}, mark=o, mark size=2pt,line
        width=1pt,#1]%
}
\newcommand{\plotlinestylethree}[1][solid]{%
    \addplot[color=blue,mark options={solid}, mark=triangle, mark size=2pt,line
        width=1pt, #1]%
}

\newcommand{\plotlinefull}[1][solid]{%
    \addplot[color=red,mark options={solid}, mark=*, mark size=2pt,line
        width=1pt,#1]%
}

\title{Matrix-Free Evaluation of High-Order Shifted Boundary Finite Element Operators}
\author{MICHAŁ Wichrowski\thanks{Interdisciplinary Center for Scientific Computing,
        Heidelberg University, Germany, \email{mt.wichrowsk@uw.edu.pl}}}

\begin{document}

\maketitle

\begin{abstract}
    This paper presents a matrix-free approach for implementing the shifted boundary method (SBM) in finite element
    analysis. The SBM is a versatile technique for solving partial differential equations on complex geometries
    by shifting boundary conditions to nearby surrogate boundaries. We focus on the efficient evaluation of shifted
    boundary operators using precomputed data and tensor-product structures. The proposed method avoids the
    explicit
    assembly of global matrices, achieving a computational complexity of $O(p^{2d-1})$ per face for the evaluation
    of
    shifted boundary contributions on elements of polynomial degree $p$ in $d$ dimensions. Numerical experiments
    validate
    the accuracy and efficiency of the approach, demonstrating its
    scalability and applicability to high-order finite element methods for both continuous and discontinuous
    Galerkin
    formulations. We compare the performance of the proposed method with a matrix-free CutFEM implementation.
\end{abstract}



\section{Introduction}

Unfitted finite element methods, which avoid body-fitted mesh generation, offer a path forward, but their efficient
implementation remains a demanding task. The Shifted Boundary Method (SBM)~\cite{main2018shifted} is an
unfitted approach that employs a structured background mesh, shifting boundary conditions from the true domain boundary
to a nearby \emph{surrogate} boundary aligned with this mesh. This strategy not only circumvents complex body-fitted
meshing  but also simplifies geometric preprocessing compared to other unfitted techniques that might involve more
intricate geometric operations.

While SBM simplifies mesh generation, its computational efficiency, particularly for high-order discretizations and
large-scale problems, can be hampered by traditional matrix-based implementations. For finite elements of polynomial
degree $p$ in $d$ spatial dimensions, the storage and manipulation of element and global matrices (with local sizes
scaling as $O(p^{2d})$) become prohibitive, often leading to memory-bound computations limited by data access rather
than arithmetic capability. This paper addresses this critical performance bottleneck by developing a matrix-free
framework specifically for SBM. The relative simplicity of SBM's geometric handling makes it a particularly promising
candidate for efficient matrix-free evaluation, aiming to unlock its full potential for high-performance computing.

Matrix-free methods provide a powerful alternative by computing the action of the finite element operator on-the-fly,
directly from the variational formulation, thus avoiding the assembly and storage of large sparse matrices. Leveraging
sum factorization techniques~\cite{kronbichler2019fast} on tensor-product basis functions common in SBM's structured
background meshes, the computational complexity for operator application can be reduced from $O(p^{2d})$ to $O(d
    p^{d+1})$ per element. This reduction is critical, as the cost of evaluating boundary terms, which scales as
$O(p^{2d-1})$, would otherwise be dominated by the $O(p^{2d})$ cost of standard matrix-based evaluation on boundary
cells. This not only drastically reduces memory requirements but also enhances arithmetic
intensity, leading to improved utilization of modern computing architectures. The primary contribution of this work is
the detailed development, analysis, and demonstration of such a matrix-free evaluation strategy for SBM operators. We
show its applicability to both continuous and discontinuous Galerkin (DG) formulations, which is an essential component
for the development of advanced, scalable iterative solvers.

The Shifted Boundary Method  provides a flexible framework for solving partial differential equations  (PDEs)  on
complex
domains without requiring body-fitted meshes. By utilizing a structured background mesh, SBM constructs a
\emph{surrogate} computational domain $\tilde{\Omega}$ composed of selected active cells, whose boundary
$\tilde{\Gamma}$ may not coincide with the true boundary $\Gamma$. Boundary conditions are imposed by transferring data
from $\Gamma$ to $\tilde{\Gamma}$, typically using closest point projections and Taylor expansions, and are enforced
weakly through Nitsche-type formulations~\cite{nitsche1971variationsprinzip, zorrilla2024shifted}.

While this avoids the complexities of generating body-fitted meshes, the primary geometric task in SBM shifts to
accurately determining the relationship between points on the surrogate boundary $\tilde{\Gamma}$ and the true boundary
$\Gamma$. The method inherently allows for the use of arbitrarily complex geometries, and crucially, avoids the need to
compute integrals over the arbitrarily shaped integration domains that arise from cell-boundary intersections. However,
for each point on $\tilde{\Gamma}$, one must find its closest point projection onto $\Gamma$, which is still a
non-trivial problem, especially for complex or implicitly defined geometries. Level set methods, which represent the
domain boundary as the zero level set of a function, are often employed in unfitted methods like SBM to facilitate
operations such as closest point projection~\cite{kuzmin2022unfitted, xue2021new}. Special treatment of
domains with corners was analyzed in~\cite{atallah2021analysis}.

Since its original formulation~\cite{main2018shifted, main2018shiftedVol2}, which restricted active cells to those
strictly inside the domain, SBM has evolved to include intersected cells based on volume fraction
criteria~\cite{yang2024optimal}. The method has been extended to high-order discretizations~\cite{atallah2022high}, and
applied to a variety of physical problems, such as Stokes flow~\cite{atallah2020analysis}, solid
mechanics~\cite{atallah2021shifted, atallah2024nonlinear}, and problems with embedded interfaces~\cite{li2020shifted,
    xu2024weighted}. In the latter, SBM has been adapted to impose jump conditions across internal boundaries,
broadening its applicability to multiphysics and multi-material scenarios. Additional developments include penalty-free
variants~\cite{collins2023penalty} and integration with level set methods for geometric
representation~\cite{kuzmin2022unfitted, xue2021new}.

Despite these advancements and the broadening scope of SBM applications, the computational cost associated with
traditional matrix-based implementations remains a significant hurdle, particularly for the high-order discretizations
and large-scale simulations where SBM's advantages are most pronounced. Addressing this efficiency challenge is
paramount to fully realize the potential of SBM. This naturally leads to the exploration of matrix-free techniques,
which form the central theme of the present work.

Matrix-free implementations are central to the efficiency of the approach presented in this work. By evaluating the
finite element operator on-the-fly, without assembling global matrices, one can achieve high performance and
scalability,
particularly for high-order discretizations and large-scale problems~\cite{kronbichler2019fast, witte2021fast}.
Matrix-free methods have been successfully applied to nonlinear problems such as
hyperelasticity~\cite{schussnig2025matrix, davydov2020matrix} and fluid-structure interaction
(FSI)~\cite{wichrowski2023exploiting}, demonstrating their potential for complex applications.
Recent advances in matrix-free techniques, such
as the use of automatic differentiation (AD) for tangent operator evaluation~\cite{wichrowski2025LargeStrain},
demonstrate that AD-based matrix-free implementations can match or even
outperform carefully tuned hand-written code. As a result, the implementation effort required for complex operators is
no longer prohibitive, making matrix-free approaches broadly accessible for advanced applications.

Applying matrix-free techniques to SBM requires careful consideration of the different terms in the variational
formulation. While standard volumetric terms on interior cells benefit directly from sum factorization, the terms
arising from the shifted boundary conditions on the surrogate boundary $\tilde{\Gamma}$ involve geometric data (closest
point projections, shift vectors) and require evaluating basis functions at points on the true boundary $\Gamma$ that
do not necessarily align with tensor-product quadrature points. Efficiently handling these boundary terms within a
matrix-free framework is crucial for the overall performance of the method.

Geometric multigrid methods~\cite{brandt1977multi, hackbusch1985multi} provide highly efficient solvers for elliptic
PDEs, and their locality makes them a natural fit for matrix-free implementations.
Recent work~\cite{wichrowski2025Geometric} has demonstrated the effectiveness of geometric multigrid solvers
specifically designed for the DG-SBM discretization, showing that such solvers can achieve mesh-independent convergence
and high parallel scalability. However, the implementation in~\cite{wichrowski2025Geometric} relies on assembled sparse
matrices. While the multiplicative Schwarz smoothers used in that work proved effective, they are challenging to
implement in a matrix-free context. To the author's disappointment, simpler additive smoothers, which are more amenable
to matrix-free implementation, were found to be insufficiently efficient for the systems arising from SBM for quadratic
and higher order elements. Therefore,
this paper focuses on the efficient matrix-free operator evaluation, which is a prerequisite for developing advanced
solvers, rather than presenting a complete multigrid solver.

Discontinuous Galerkin (DG) methods~\cite{reed1973triangular, zienkiewicz2003discontinuous, arnold1982interior,
    cockburn2000development} employ discontinuous polynomial basis functions, enabling element-wise independence and
local conservation. For SBM, DG discretizations are attractive due to their natural handling of weak continuity and
boundary conditions, and their suitability for efficient, cell-local multigrid smoothers. The geometric multigrid
solver in~\cite{wichrowski2025Geometric} exploits these features for robust and scalable DG-SBM solvers. While DG
methods introduce more degrees of freedom than continuous Galerkin approaches, their flexibility and compatibility
with matrix-free multigrid make them highly effective for high-performance SBM implementations.

Among unfitted finite element methods, CutFEM~\cite{burman2015cutfem} stands out as the closest competitor to SBM,
offering a general framework for discretizing PDEs on complex domains by directly cutting the background mesh to fit
the physical geometry. CutFEM operates by integrating over the true domain, requiring specialized quadrature rules for
intersected (cut) cells and faces. This approach enables high geometric flexibility but introduces challenges such as
handling small cut elements, which can lead to ill-conditioning. To address this, stabilization techniques like the
ghost penalty~\cite{burman2010ghost, wichrowski2025matrix} are essential, though they may introduce issues such as
locking if not carefully designed~\cite{badia2022linking, bergbauer2024high, burman2022design}.

CutFEM has been successfully applied to a wide range of problems, including Stokes flow~\cite{burman2014fictitious},
elasticity~\cite{hansbo2017cut}, and two-phase flows~\cite{claus2019cutfem}. Discontinuous Galerkin (DG) formulations
have also been combined with CutFEM~\cite{gurkan2019stabilized, bergbauer2024high}, further enhancing its flexibility.

Matrix-free implementations of CutFEM have recently been developed~\cite{bergbauer2024high, wichrowski2025matrix},
enabling efficient operator application even for high-order elements and large-scale problems. The matrix-free
evaluation of CutFEM was first described in~\cite{bergbauer2024high}, focusing on efficient computation of cut cell
contributions. In~\cite{wichrowski2025matrix}, the method was extended by
implementing matrix-free evaluation of the ghost penalty based on tensor products, improving the
computational efficiency of the stabilization terms. However, the irregular integration domains and the need for
complex geometric queries in cut cells make matrix-free CutFEM more challenging and potentially less efficient than

The choice between SBM and body-fitted methods often depends on the application. For problems involving moving or
evolving boundaries, such as in fluid-structure interaction or shape optimization, the cost of repeatedly generating
high-quality body-fitted meshes can be prohibitive. In these scenarios, SBM's ability to handle complex geometries on a
fixed background mesh offers a decisive advantage, simplifying the overall simulation workflow.

Regarding preconditioning, CutFEM presents additional difficulties. While optimal preconditioners have been
proposed~\cite{gross2023analysis, gross2021optimal}, achieving mesh-independent convergence, iteration counts can
remain high. In~\cite{bergbauer2024high}, a multigrid preconditioner based on a cell-wise Additive Schwarz smoother was
used for DG-CutFEM, showing promise in a matrix-free context. Nevertheless, the smoothing step typically requires a
relatively large number of matrix-vector products, impacting overall efficiency.

This paper presents a detailed matrix-free framework for evaluating the finite element operators arising from SBM
discretizations, applicable to formulations with both continuous and discontinuous elements. We discuss
the data structures required to store geometric information for the shifted boundary terms and detail how sum
factorization is applied to different parts of the operator evaluation (interior cells, interior faces for DG, and
surrogate boundary faces). We analyze the computational complexity of each component and demonstrate how this
matrix-free approach enables the efficient application of SBM operators, paving the way for scalable solvers like
geometric multigrid. We also provide a conceptual comparison of the local computational cost of matrix-free SBM
evaluations with those of matrix-free CutFEM, highlighting the potential advantages of SBM due to its simpler
integration domains.

The remainder of this paper is organized as follows. Section~\ref{sec:sbm_formulation} introduces the SBM formulation
and its weak form for both continuous and discontinuous Galerkin discretizations. Section~\ref{sec:matrix_free_eval}
details the matrix-free implementation, covering data structures, evaluation strategies for different operator
components, and the exploitation of tensor-product structures through sum factorization. In
Section~\ref{sec:cutfem_introduction}, we provide a brief overview of CutFEM. Section~\ref{sec:numerical_results}
presents comprehensive numerical experiments that validate the method's efficiency, including microbenchmarks that
measure local computational costs for individual mesh entities, parallel scalability tests, and performance comparisons
with matrix-free CutFEM implementations. Finally, Section~\ref{sec:conclusion} summarizes the key findings and outlines
directions for future research.

\section{Shifted Boundary Method Formulation}\label{sec:sbm_formulation}

\begin{figure}[!ht]
    \centering
    \begin{tikzpicture}[scale=1.7]
    \draw[step=1cm, gray, very thin] (0.8,0) grid (5.2,2.2);

    \foreach \x/\y in {2/1, 1/0, 1/1, 2/0, 3/0, 4/0}
        {
            \fill[green, opacity=0.5] (\x,\y) rectangle (\x+1,\y+1);
        }
    \foreach \x/\y in {1/0 }{
            \fill[green, opacity=0.5] (\x,\y) rectangle (\x-0.2,\y+2);
        }

    \draw[line width=2pt, red, domain=51:90, samples=100, variable=\t]
    plot ({0.8 + 7.0*cos(\t)}, {-4.4 + 7.0*sin(\t)});
    \node at (2.6,2.5) {$\Gamma$};

    \fill[green!30, opacity=0.4] plot[domain=51:90, samples=100, variable=\t] ({0.8 + 7.0*cos(\t)}, {-4.4 +
            7.0*sin(\t)}) -- (0.8,0) -- (5.2,0) -- cycle;

    \pgfmathsetmacro{\myangle}{66}
    \draw[<-, thick] ({0.8 + 6.9*cos(\myangle) + 0.3}, {-4.4 + 6.9*sin(\myangle)}) -- ({0.8 + 5.9*cos(\myangle)
            +
            0.3}, {-4.4 +
            5.9*sin(\myangle)});
    \pgfmathsetmacro{\labelangle}{\myangle - 6} 
    \node at ({0.8 + 6.9*cos(\labelangle) - 0.6}, {-4.4 + 6.9*sin(\labelangle)+0.1}) {\large $\mathbf{d}$};

    \draw[->, thick] ({0.8 + 6.8*cos(\myangle-5) + 0.3}, {-4.4 + 6.8*sin(\myangle-5)}) -- ({0.8 +
            7.6*cos(\myangle-5)
            +
            0.3}, {-4.4 +
            7.6*sin(\myangle-5)});
    \node at ({0.8 + 7.5*cos(\labelangle-5) - 0.54}, {-4.4 + 7.5*sin(\labelangle-5)+0.36}) {\large $\mathbf{n}$};

    \draw[<-, thick]  ({0.8 + 5.9*cos(\myangle) +
            0.3}, {-4.4 + 6.9*sin(\myangle)}) -- ({0.8 + 5.9*cos(\myangle) +
            0.3}, {-4.4 +
            5.9*sin(\myangle)});
    \node at (3.4, 1.5) {\large $\tilde{\mathbf{n}}$}; 

    \draw[line width=2pt, blue] (0.8,2) -- (3,2) -- (3,1) -- (4,1) -- (5,1)-- (5,0);

    \draw[line width=2pt, blue, dashed] (3,2) -- (4,2) -- (4,2) -- (4,1) ;
    \node at (1.7,2.2) {$\tilde{\Gamma}$};

    \node at (2.5,0.7) {$\tilde{\Omega}$};

    \node at (1.2,2.3) {${\Omega}$};

\end{tikzpicture}
    \caption{Schematic illustrating the background mesh, interior cells (green), the surrogate boundary
        $\tilde{\Gamma}$ (thick blue line) along the upper boundary of the interior cells, and the true boundary
        $\partial\Omega$ (red).\label{fig:sbm_setup}}
\end{figure}
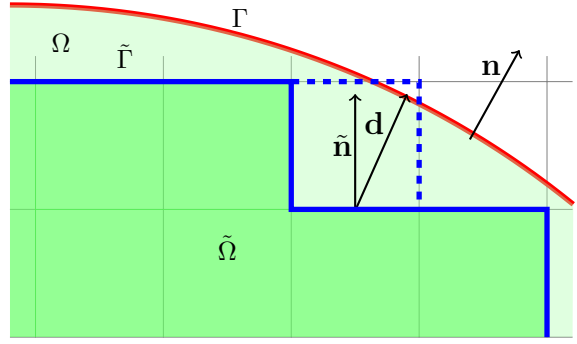

Consider a bounded Lipschitz domain $\Omega \subset \mathbb{R}^d$ with boundary $\partial\Omega = \Gamma$. We aim
to solve the model Poisson problem:
\[
    -\Delta u = f \quad \text{in } \Omega, \quad u = g \quad \text{on } \Gamma.
\]
The Shifted Boundary Method (SBM) addresses this problem by replacing the original physical domain $\Omega$ with a
surrogate computational domain $\tilde{\Omega}$. This surrogate domain (illustrated in Figure~\ref{fig:sbm_setup}) is
constructed as a union of cells from a fixed background mesh $\mathcal{T}_h$ that does not necessarily conform to the
true boundary $\Gamma$. This approach circumvents the need for body-fitted mesh generation, which can be complex and
time-consuming, especially for intricate or evolving geometries \cite{main2018shifted}.
Typically, $\tilde{\Omega}$ comprises cells from $\mathcal{T}_h$ designated as \emph{active}. Active cells might be
those lying entirely within $\Omega$ or those intersecting $\Omega$ significantly, based on a chosen criterion (e.g.,
volume fraction \cite{yang2024optimal}). The boundary of this surrogate domain is denoted by $\tilde{\Gamma} =
    \partial\tilde{\Omega}$.

The  SBM formulation transfers the boundary conditions from the true boundary $\Gamma$ to the surrogate boundary
$\tilde{\Gamma}$. For each point $\boldsymbol{x}_s \in \tilde{\Gamma}$, a corresponding point $\boldsymbol{x}$
on the true boundary $\Gamma$ is identified (commonly via closest point projection). The boundary condition
$g(\boldsymbol{x})$ at this true boundary point is then extrapolated or transferred to $\boldsymbol{x}_s$ to define a
surrogate boundary condition $g_s(\boldsymbol{x}_s)$, often using techniques like Taylor series expansions.

The problem is then formulated weakly on the surrogate domain $\tilde{\Omega}$. Using a suitable finite element space
$V_h \subset H^1(\tilde{\Omega})$, the surrogate boundary conditions are enforced using a Nitsche-type
method~\cite{nitsche1971variationsprinzip}. The SBM seeks $u_h \in V_h$ such that for all test functions $v_h \in V_h$:
\[
    \int_{\tilde{\Omega}} \nabla u_h \cdot \nabla v_h \, d\boldsymbol{x} - \int_{\tilde{\Gamma}}
    \frac{\partial u_h}{\partial n_s} v_h \, ds - \int_{\tilde{\Gamma}} (u_h - g_s) \frac{\partial v_h}{\partial
        n_s} \, ds + \int_{\tilde{\Gamma}} \frac{\beta}{h_s} (u_h - g_s) v_h \, ds = \int_{\tilde{\Omega}} f v_h \,
    d\boldsymbol{x}.
\]
Here, $n_s$ represents the outward unit normal vector to the surrogate boundary $\tilde{\Gamma}$, $h_s$ is a local
mesh size parameter associated with $\tilde{\Gamma}$, and $\beta > 0$ is a penalty parameter that must be chosen
sufficiently large to ensure stability of the formulation.

The specific Nitsche terms used to enforce the boundary condition can vary. The formulation presented above is one
common choice. Another widely used variant, known for its symmetry, is given by:
\[
    \int_{\tilde{\Omega}} \nabla u_h \cdot \nabla v_h \, d\boldsymbol{x} + \int_{\tilde{\Gamma}} \left(
    \sigma (u_h - g_s) v_h - (u_h - g_s) \frac{\partial v_h}{\partial n_s} - \frac{\partial u_h}{\partial n_s} v_h
    \right) \, ds = \int_{\tilde{\Omega}} f v_h \, d\boldsymbol{x},
\]
where $\sigma$ is a penalty parameter, typically chosen to be sufficiently large (e.g., $\sigma \approx C/h_s$ for
some constant $C > 0$) to ensure stability.

\subsection{Shifting boundary conditions}\label{ssec:shifting_boundary_conditions}
The Dirichlet condition $g$ on the true boundary $\Gamma$ is transferred to the surrogate boundary $\tilde{\Gamma}$ to
define $g_s$. This process involves an \emph{extension operator} $\mathcal{E}$ that extrapolates boundary conditions
from $\Gamma$ to $\tilde{\Gamma}$.

We assume that the Dirichlet boundary condition $g$ is given as a restriction of a function $u^\star$ defined on the
entire domain $\Omega$ to the boundary $\Gamma$. While the choice of this function $u^\star$ (as an extension of $g$
from $\Gamma$ into $\Omega$) is not unique, we take $u^\star$ to be equal to the solution $u$ in the surrogate domain
$\tilde{\Omega}$.

For each point $\tilde{\mathbf{x}} \in \tilde{\Gamma}$, let $\mathbf{x} \in \Gamma$ be its closest point projection
onto the true boundary, and let $\mathbf{d} = \mathbf{x} - \tilde{\mathbf{x}}$ be the shift vector. The function $u$ is
extended from $\Gamma$ to $\tilde{\Gamma}$ using a Taylor expansion:
\[
    \mathcal{E}u^\star(\tilde{\mathbf{x}}) = u^\star(\mathbf{x}) - \mathbf{d} \cdot \nabla u (\tilde{\mathbf{x}}) +
    \cdots
\]
where $\mathcal{E}u^\star$ denotes the extrapolated boundary condition. The function $u^\star$ is assumed to be smooth
in a neighborhood of $\tilde{\Gamma}$, which allows for the Taylor expansion to be valid.
By substituting this into the weak formulation, we obtain a variational formulation for the shifted boundary problem
with the extrapolated boundary condition on $\tilde{\Gamma}$ enforced in a Nitsche-like manner.

The SBM weak formulation seeks $u_h \in V_h$ such that for all $v_h \in V_h$:
\begin{equation}
    \begin{split}
        \int_{\tilde{\Omega}} \nabla u_h \cdot \nabla v_h \, dx
        - \int_{\tilde{\Gamma}} (\nabla u_h \cdot \tilde{\mathbf{n}}) v_h \, ds
        -  \int_{\tilde{\Gamma}} (\nabla v_h \cdot \tilde{\mathbf{n}}) \; \mathcal{E}u_h \, ds
        + \int_{\tilde{\Gamma}} \sigma_\Gamma \; \mathcal{E} u_h  \; v_h \, ds = \\
        = \int_{\tilde{\Omega}} f v_h \, dx
        - \int_{\tilde{\Gamma}} (\nabla v_h \cdot \tilde{\mathbf{n}}) g \, ds
        + \int_{\tilde{\Gamma}} \sigma_\Gamma g v_h \, ds,
    \end{split}
\end{equation}
where $\tilde{\mathbf{n}}$ is the outward normal to $\tilde{\Gamma}$ and $\sigma_\Gamma$ is a penalty parameter.

In the matrix-free implementation, the extension operator $\mathcal{E}$ is applied directly to the discrete solution
$u_h$. Since $u_h$ is a piecewise polynomial function defined on the background mesh $\mathcal{T}_h$, its Taylor
expansion can be computed directly by evaluating the function values and gradients at points on $\tilde{\Gamma}$ and
applying the shift $\mathbf{d}$. This allows efficient evaluation of the extrapolated boundary values $\mathcal{E}u_h$
during the matrix-free operator application without requiring higher-order derivatives. The geometric data that must be
precomputed and stored includes the locations of corresponding points on the true boundary $\Gamma$, stored in
reference coordinates for efficient lookup. For each quadrature point on the surrogate boundary face, we store the
reference coordinates of its closest point projection onto $\Gamma$, along with the shift vector $\mathbf{d}$ and any
required boundary data values.

While the above formulation is presented in the context of continuous finite element spaces $V_h \subset
    H^1(\tilde{\Omega})$, the SBM can also be effectively discretized using Discontinuous Galerkin (DG)
methods~\cite{wichrowski2025Geometric}. DG methods
employ basis functions that are piecewise polynomials, discontinuous across element interfaces. This inherent
discontinuity offers greater flexibility, particularly in handling complex geometries and designing robust smoothers
for multigrid solvers. In the SBM-DG context, the Nitsche-type boundary condition enforcement on $\tilde{\Gamma}$
remains similar, but additional terms arise from penalizing jumps across interior faces $\mathcal{F}$ of
$\mathcal{T}_h$ within
$\tilde{\Omega}$.  Namely, the standard Laplacian term is replaced with the DG formulation that includes interior
penalty terms. The weak form becomes:
\[
    \int_{\tilde{\Omega}} \nabla u_h \cdot \nabla v_h \, dx \rightarrow \sum_{K \in \mathcal{T}_h} \int_{K} \nabla u_h
    \cdot \nabla v_h \, dx + \sum_{F \in \mathcal{F}} \int_{F} \left( \sigma_F [u_h] \cdot [v_h] - \{ \nabla u_h \}
    \cdot
    [v_h] - [u_h] \cdot \{ \nabla v_h \} \right) \, ds,
\]
where $[u_h]$ denotes the jump of $u_h$ across the interior face $F$, $\{\nabla u_h\}$ denotes the average of the
gradient across $F$, $\mathcal{F}$ is the set of interior faces, and $\sigma_F$ is the penalty parameter.

The element-local nature of DG discretizations is particularly advantageous for constructing
efficient cell-wise smoothers within a geometric multigrid framework, which is crucial for tackling the ill-conditioned
systems often produced by SBM.
\section{Matrix-Free Evaluation of SBM Operators}\label{sec:matrix_free_eval}

Building upon the SBM formulation presented in the previous section, we now detail the matrix-free approach for
evaluating the resulting finite element operators. The action of the operator on a vector, such as a trial solution
vector in an iterative solver, is computed on-the-fly by summing contributions from different parts of the
computational domain $\tilde{\Omega}$. This evaluation process distinguishes between contributions from interior cells
(those within $\tilde{\Omega}$ not adjacent to the surrogate boundary $\tilde{\Gamma}$), interior faces (for
Discontinuous Galerkin discretizations, these are faces shared by two cells within $\tilde{\Omega}$), and surrogate
boundary faces (faces of cells in $\tilde{\Omega}$ that lie on $\tilde{\Gamma}$, where the shifted boundary conditions
are applied).

The evaluation process can be decomposed into distinct computational components, each requiring specific treatment
within the matrix-free framework. We first examine the evaluation strategies for different types of terms in the SBM
formulation, considering how tensor-product structures can be leveraged for computational efficiency. Subsequently, we
discuss the data structures and storage requirements necessary to support these evaluation strategies, particularly
focusing on the geometric information needed for the shifted boundary condition terms and the precomputed data that
enables efficient on-the-fly computation.

\subsection{Tensor-Product Structure and Sum Factorization}\label{ssec:tensor_product}

\begin{figure}[!ht]
    \centering
    \begin{tikzpicture}[scale=1.5]
        \draw[line width=1pt, dashed] (-2,0) rectangle (0,2) node[midway, below] {};
        \foreach \x in {0,1,2}
        \foreach \y in {0,1,2 }
        {
        \node[circle,inner sep=2pt,fill=red!70,label={[label distance=2pt]below left:{\small (\x,\y)}}]
        (K1-\x-\y) at (-2+\x*2/2,\y*2/2) {};
        \pgfmathtruncatemacro{\index}{\x + \y * 5}
        \node[above right=1pt,font=\small] at (K1-\x-\y) {\index};
        }

    \end{tikzpicture}
    \caption{Tensor product numbering of degrees of freedom for a quadratic element $K$ in 2D.
        The numbers indicate the lexicographical ordering of the DoFs within the cell, and the pairs in parentheses
        denote the corresponding multi-indices $(i_1, i_2)$ for $p=2$.
        This ordering is fundamental to the efficiency of sum factorization, as it allows multi-dimensional operations
        to be decomposed into a sequence of one-dimensional sweeps.
        \label{fig:tensor_product_numbering}
    }
\end{figure}
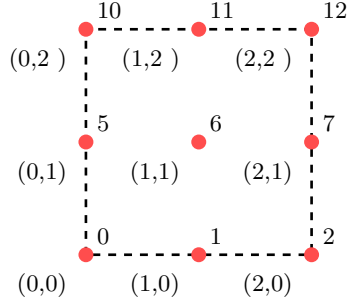

On Cartesian elements, which are natural for the background mesh in SBM, the basis functions are typically constructed
as tensor products of one-dimensional polynomials. Let $K$ be a $d$-dimensional Cartesian cell. The degrees of freedom
within this cell are numbered lexicographically, as illustrated for $d=2$ in
Figure~\ref{fig:tensor_product_numbering}. A multi-index $\boldsymbol{i} = (i_1, \dots, i_d)$, where $0 \le i_\ell \le
    p$ for polynomial degree $p$, can be used to identify each basis function
$\psi_{\boldsymbol{i}}(\boldsymbol{x})$.
These functions are formed by the product of $d$ one-dimensional basis functions $\hat{\phi}_{p}(\xi_l)$:
\[
    \psi_{\boldsymbol{i}}(\boldsymbol{x}) = \prod_{l=1}^d \hat{\phi}_{i_l}(\xi_l(\boldsymbol{x})),
\]
where $\boldsymbol{\xi}(\boldsymbol{x})$ maps physical coordinates $\boldsymbol{x} \in K$ to reference coordinates
$\boldsymbol{\xi} \in \hat{K} = {[0,1]}^d$.

The evaluation of terms like $\int_K \nabla u_h \cdot \nabla v_h , d\boldsymbol{x}$ involves three key operations:
first, computing values (e.g., gradients) of the function $u_h$ from its coefficients at all quadrature points; second,
performing pointwise operations within a loop over these
quadrature points; and third, integrating these results by summing over quadrature points and testing against basis
functions. A naive approach to the first operation, evaluating,
for instance, the gradient of $u_h = \sum_{\boldsymbol{j}} u_{\boldsymbol{j}} \psi_{\boldsymbol{j}}$ at $N_q$
quadrature points would involve summing over all $(p+1)^d$ basis functions for each quadrature point, leading to a
cost of $O(N_q (p+1)^d)$. If $N_q$ is chosen as $O((p+1)^d)$, a natural choice being $(p+1)^d$ Gauss points for exact
integration of certain terms, this becomes $O((p+1)^{2d})$.

However, by choosing a quadrature rule that itself has a tensor-product structure (e.g., using $(p+1)$ Gauss points
in each of the $d$ directions, totaling $N_q = (p+1)^d$ points), sum factorization can be employed to dramatically
reduce the cost~\cite{kronbichler2012generic}. Sum factorization breaks down the multi-dimensional evaluation into a
sequence of $d$ one-dimensional operations. For example, to evaluate $u_h$ and its derivatives at all $(p+1)^d$
tensor-product quadrature points, one applies 1D evaluation/differentiation operators along each coordinate direction
sequentially. This reduces the cost of obtaining all gradient components at all quadrature points to $O(d
    (p+1)^{d+1})$.

Once the gradients of $u_h$ are available at the quadrature points, the subsequent integration step, which
involves contracting these with test function gradients (for the stiffness matrix term) and summing over quadrature
points, can also be structured to exploit sum factorization. The operation of testing and integrating (e.g., computing
$\int_K \nabla \psi_{\boldsymbol{i}} \cdot
    \nabla u_h \, d\boldsymbol{x}$ for all $\boldsymbol{i}$) is structurally similar to the evaluation step (transposed
evaluation). Thus, this "test and integrate" phase also achieves a complexity of $O(d (p+1)^{d+1})$. The procedure is
summarized in Algorithm~\ref{alg:mf_cell}.
\begin{algorithm2e}[!ht]
    \SetKwInOut{Input}{Given}
    \SetKwInOut{Output}{Return}
    \Input{${u}$ -- current FE solution}
    \Output{${w}=	\mathcal{A}_K u$}

    Gather element-local vector values from cell $K$\;
    Evaluate $\nabla u_h$ at all quadrature points $q$ on $K$: \\
    \Indp
    $\{\nabla u_h(x_q)\}_{q \in K}$ \tcp*{Sum factorization}
    \Indm
    \ForEach{quadrature point $q$ on $K$ }{
        Compute pointwise integrand contribution $I_q = \nabla u_h(x_q) $\;
        Submit $I_q$ for integration against test function gradients\;
    }
    Integrate submitted contributions: \\
    \Indp
    $w_{Ki} \leftarrow \sum_q \nabla \phi_i(x_q) \cdot I_q$ for each basis function $\phi_i$ \tcp*{Sum factorization}
    \Indm
    Scatter results to ${w}$

    \caption{Local evaluation of cell contributions to the SBM operator. \label{alg:mf_cell}}
\end{algorithm2e}

\subsection{Evaluation of Interior Face Terms (DG)}\label{ssec:dg_faces}
In Discontinuous Galerkin (DG) formulations, additional terms arise from integrals over interior faces $F$ (i.e., faces
shared by two cells within $\tilde{\Omega}$, so $F \not\subset \tilde{\Gamma}$). These terms are essential for weakly
enforcing continuity or penalizing discontinuities. A common form for these interior face terms, for instance in the
Symmetric Interior Penalty Galerkin (SIPG) method, is:
\[
    \sum_{F \in \mathcal{F}_{\text{int}}} \int_F \left( \sigma_F \llbracket u_h \rrbracket \cdot \llbracket v_h
    \rrbracket - \{ \nabla u_h \} \cdot \llbracket v_h \rrbracket - \llbracket u_h \rrbracket \cdot \{ \nabla v_h \}
    \right) \, ds,
\]
where $\mathcal{F}_{\text{int}}$ is the set of interior faces, $\llbracket \cdot \rrbracket$ denotes the jump operator
(e.g., $u_{h,K_1} - u_{h,K_2}$ across the face between cells $K_1$ and $K_2$), $\{ \cdot \}$ denotes the average
operator (e.g., $0.5(\nabla u_{h,K_1} + \nabla u_{h,K_2})$), and $\sigma_F$ is a penalty parameter.

The matrix-free evaluation of these interior face terms, detailed in Algorithm~\ref{alg:mf_face}, also relies heavily
on sum factorization for its efficiency. For each interior face $F$ shared by cells $K_1$ and $K_2$, the process starts
by evaluating the traces (values and gradients) of the solution $u_h$ from both cells at all quadrature points on $F$,
using sum factorization to efficiently compute $u_{h,K_1}$, $\nabla u_{h,K_1}$, $u_{h,K_2}$, and $\nabla u_{h,K_2}$ at
these points from the local coefficient vectors.

With these values available at each quadrature point on the face, the algorithm proceeds to compute the necessary jump
$\llbracket u_h \rrbracket(x_q) = u_{h, K_1}(x_q) - u_{h, K_2}(x_q)$ and the average of the normal derivative $\{
    \nabla u_h \}(x_q) \cdot \mathbf{n} = 0.5 (\nabla u_{h, K_1}(x_q) + \nabla u_{h, K_2}(x_q)) \cdot \mathbf{n}$.
These
quantities, along with the face penalty parameter $\sigma_F$, are then used to construct the pointwise contributions to
the DG face integrals that will be tested against the basis functions. For example, for the SIPG terms involving
testing against the value of the test function $v_h$, a term like $T_v = \sigma_F \llbracket u_h \rrbracket(x_q) - \{
    \nabla u_h \}(x_q) \cdot \mathbf{n}$ is formed. For terms involving testing against the normal derivative of $v_h$,
a
term proportional to $\llbracket u_h \rrbracket(x_q)$ is formed. These pointwise integrands are then submitted for
integration.

Finally, these submitted pointwise contributions are integrated against the corresponding traces of the test functions
$v_h$ (i.e., their values and normal derivatives from cells $K_1$ and $K_2$ on the face $F$). This integration step is
again performed efficiently using sum factorization, effectively applying a transposed evaluation operation. The
resulting local force vector contributions for cells $K_1$ and $K_2$ are then scattered to the global output vector
$w$.

Since interior faces are $(d-1)$-dimensional and typically aligned with coordinate axes in structured Cartesian meshes,
they inherit tensor-product structure. The evaluation of traces from cell data to face quadrature points and the
subsequent integration steps (testing against basis functions) leverage sum factorization over these $(d-1)$
dimensions. The overall computational complexity for processing one interior face, considering the operations related
to the $(p+1)^d$ degrees of freedom of the two adjacent cells, is approximately $O(d (p+1)^d)$. This is because the
evaluation of values/gradients from the $d$-dimensional cell data to the $(d-1)$-dimensional face quadrature points,
and the corresponding integration step, are the dominant costs.

\begin{algorithm2e}[!ht]
    \SetKwInOut{Input}{Given}
    \SetKwInOut{Output}{Return}
    \Input{  ${u}$ -- current FE solution \\}
    \Output{ ${w}=	\mathcal{A}_F u$ }
    Gather element-local vector values from cell $K_1$ and cell $K_2$ adjacent to the face $F$\\
    Evaluate traces of functions and gradients at face quadrature points for both cells: \\
    \Indp
    $\{u_{i,K_1}(x_q), \nabla u_{i,K_1}(x_q)\}_{q \in F, i \in K_1}$ \tcp*{Sum factorization}
    $\{u_{j,K_2}(x_q), \nabla u_{j,K_2}(x_q)\}_{q \in F, j \in K_2}$ \tcp*{Sum factorization}
    \Indm
    Compute penalty parameter $\sigma_F$\\
    \ForEach{quadrature point $q$ on $F$ }{
        Compute jump $\llbracket u_h \rrbracket(x_q) = u_{h, K_1}(x_q) - u_{h, K_2}(x_q)$\;
        Compute average normal derivative $\{ \nabla u_h \}(x_q) \cdot \mathbf{n} = 0.5 (\nabla u_{h,
                    K_1}(x_q) + \nabla u_{h, K_2}(x_q)) \cdot \mathbf{n}$\;
        Compute value contribution term $T_v = \sigma_F \llbracket u_h \rrbracket(x_q) - \{ \nabla u_h \}(x_q)
            \cdot \mathbf{n}$\;
        Submit $T_v$ to cell $K_1$'s value contribution\;
        Submit $-T_v$ to cell $K_2$'s value contribution\;
        Submit $-0.5 \llbracket u_h \rrbracket(x_q)$ to cell $K_1$'s normal derivative contribution\;
        Submit $-0.5 \llbracket u_h \rrbracket(x_q)$ to cell $K_2$'s normal derivative contribution\;
    }
    Integrate contributions  of cell  $K_1$  \tcp*{Sum factorization}
    Integrate contributions  of cell $K_2$	\tcp*{Sum factorization}
    Scatter results to ${w}$

    \caption{Local evaluation of internal face contributions to the DG-SBM operator. \label{alg:mf_face}}
\end{algorithm2e}

\subsection{Evaluation of Shifted Boundary Condition Terms}\label{ssec:sbm_bc_terms}
The evaluation of shifted boundary condition terms, detailed in Algorithm~\ref{alg:mf_sbm_face}, combines the
tensor-product structure of basis functions on the surrogate boundary with precomputed geometric data to handle the
shift. This algorithm outlines
this process, which, while analogous to the treatment of interior face terms in DG methods, includes distinct steps to
accommodate the geometry of the shifted boundary.

The algorithm begins by gathering the local degrees of freedom from the cell $K$ adjacent to $F_s$. Quantities like
$\nabla u_h$ required directly at the quadrature points $\tilde{\mathbf{x}}_q$ on $F_s$ (e.g., for terms like
$\int_{F_s} (\nabla u_h \cdot \tilde{\mathbf{n}}) v_h \, ds$) are evaluated efficiently using sum factorization over
the $(d-1)$ dimensions of the face.

The evaluation of the extension operator $\mathcal{E}u_h$ at each quadrature point $\tilde{\mathbf{x}}_q$ on $F_s$
requires
special attention. In this context, $\mathcal{E}u_h(\tilde{\mathbf{x}}_q)$ is taken as the value of the solution $u_h$
at a corresponding shifted point $\mathbf{x}_q$ on  the true boundary $\Gamma$. The precomputed
reference coordinates of $\mathbf{x}_q$ (using coordinates in the reference cell $K$)  are retrieved for
each $\tilde{\mathbf{x}}_q$. It is important to note that the point $\mathbf{x}_q$ is outside of the cell $K$.

Since the locations of these shifted points $\mathbf{x}_q$ do not generally form a
tensor-product structure within their respective cells, the evaluation of $u_h(\mathbf{x}_q)$ cannot leverage sum
factorization for this specific step. Instead, for each $\tilde{\mathbf{x}}_q$, the value $u_h(\mathbf{x}_q)$ is
computed by a standard point evaluation: summing contributions from all $(p+1)^d$ basis functions of the cell
containing $\mathbf{x}_q$. This operation has a computational cost of $O((p+1)^d)$ per point $\mathbf{x}_q$.
Given that there are typically $O((p+1)^{d-1})$ quadrature points on the $(d-1)$-dimensional face $F_s$, the total
complexity for evaluating $\mathcal{E}u_h$ at all quadrature points across one such surrogate boundary face becomes
$O((p+1)^{d-1} \cdot (p+1)^d) = O((p+1)^{2d-1})$. This cost is notably higher than the $O(d(p+1)^d)$ complexity for
evaluating $u_h$ on $F_s$ using sum factorization (to get values from cell $K$ onto its face $F_s$) and can become a
dominant factor for high polynomial degrees $p$ or in
higher dimensions $d$.

The integrand contributions for the SBM boundary terms are then assembled at each quadrature point. These typically
include terms involving the normal derivative of $u_h$, the penalty term proportional to $\mathcal{E}u_h$, and possibly
other Nitsche-type contributions, all evaluated at the shifted boundary location. Each contribution is submitted for
integration against the appropriate test function traces (values or normal derivatives), again using sum factorization
for efficiency.

After integration, contributions are scattered to the global output vector. The dominant cost for surrogate boundary
faces is evaluating $u_h$ at the $O((p+1)^{d-1})$ shifted points $\mathbf{x}_q$ to compute $\mathcal{E}u_h$. Other
operations on $F_s$, such as evaluating $\nabla u_h$ and integrating submitted terms against test functions, leverage
sum factorization over the $(d-1)$ dimensions of the face and have a complexity of approximately $O(d(p+1)^d)$.
Consequently, the overall workload for surrogate boundary faces, $O((p+1)^{2d-1})$, is typically higher
than that for interior faces (which cost $O(d(p+1)^d)$), particularly for $d=3$ or for high polynomial degrees $p$.

\begin{algorithm2e}[!ht]
    \SetKwInOut{Input}{Given}
    \SetKwInOut{Output}{Return}
    \Input{  ${u}$ - current FE solution }
    \Output{ ${w}=	\mathcal{A}_{F_s} u$ }
    Gather element-local vector values from the cell $K$ adjacent to the surrogate boundary face $F_s$\;

    Evaluate $\nabla u_h$ at face quadrature points $\mathbf{x}_q$  \tcp*{Sum factorization}
    \ForEach{quadrature point $\tilde{\mathbf{x}}_q$ on $F_s$ }{
        Retrieve precomputed reference coordinates of the shifted point $\mathbf{x}_q$ on $\Gamma$ corresponding to
        $\tilde{\mathbf{x}}_q$ on $F_s$\;
        Evaluate $u_h$ at the shifted point $\mathbf{x}_q$ on $\Gamma$\;
        Compute value integrand contribution $I_{value} = - (\nabla u_h(\tilde{\mathbf{x}}_q) \cdot
            \tilde{\mathbf{n}}_q) +
            \sigma_\Gamma \; \mathcal{E} u_h(\mathbf{x}_q)$\;
        Compute gradient integrand contribution $I_{grad} = - \mathcal{E}u_h(\mathbf{x}_q) $\;
        Submit $I_{value}$ for integration against $v_h$\;
        Submit $I_{grad}$ for integration against $\nabla v_h \cdot \tilde{\mathbf{n}}$\;
    }
    Integrate submitted contributions \tcp*{Sum factorization}
    Scatter results into ${w}$\;

    \caption{Local evaluation of surrogate boundary face contributions to the SBM operator.}
    \label{alg:mf_sbm_face}
\end{algorithm2e}

\subsection{Data Structures for Matrix-Free SBM}\label{ssec:data_structures}
For cells $K \in \tilde{\Omega}$ that are located in the interior of the surrogate domain, meaning they are not
adjacent to the surrogate boundary $\tilde{\Gamma}$, significant optimizations are possible. If the background mesh is
uniform, these interior cells are often identical up to translation and scaling. Consequently, geometric information
such as the Jacobians of the mapping from a reference cell, as well as the quadrature rules, can be precomputed once
and then reused for all such standard interior cells.

Faces $F_s$ on the surrogate boundary $\tilde{\Gamma}$ are treated differently because they need specific geometric
data to handle the shifted boundary conditions. For each quadrature point $\tilde{\mathbf{x}}$
on such a face $F_s \subset \tilde{\Gamma}$, we precompute and store information about
its corresponding
point $\mathbf{x}$ on the true physical boundary $\Gamma$. The evaluation of the SBM boundary terms
typically involves the value of the solution at this true boundary point $\mathbf{x}$.
For this evaluation, and for computing the shape functions at $\mathbf{x}$, the coordinates of
this corresponding point $\mathbf{x}$ are stored in the reference coordinate system of the cell on $\Gamma$ that
contains $\mathbf{x}$.

\subsection{Matrix-Free Operator Application Workflow}\label{ssec:mf_operator_workflow}
The matrix-free operator application for the SBM system proceeds by traversing all relevant mesh entities (cells and
faces) and invoking specialized evaluation kernels tailored to each entity type and its role in the discretization. The
process is structured to maximize computational efficiency and to facilitate parallel execution, particularly in
distributed-memory environments using MPI with non-blocking communication.

The algorithm begins by initializing the global output vector $w$ to zero. Then a loop over all cells in
$\tilde{\Omega}$ computes the volumetric contributions (Algorithm~\ref{alg:mf_cell}). This is followed by a loop over
all faces. Inside this face loop, a distinction  is made between interior faces (for DG formulations, see
Algorithm~\ref{alg:mf_face}) and surrogate boundary faces (for SBM boundary conditions, see
Algorithm~\ref{alg:mf_sbm_face}), and the appropriate kernel is called. After all local computations are complete, the
necessary data is exchanged between MPI processes to finalize the global result vector.

The computational workload is not uniform across all cells in the surrogate domain $\tilde{\Omega}$; cells adjacent to
the surrogate boundary $\tilde{\Gamma}$ incur a significantly higher computational cost due to the evaluation of the
shifted boundary terms (Algorithm~\ref{alg:mf_sbm_face}). To distribute the workload evenly among MPI processes, we
assign a weight to each cell that reflects its computational cost. Cells outside the computational domain are assigned
a weight of zero. Interior cells are given a baseline weight equal to 10 corresponding to the cost of volumetric
integration. For cells adjacent to the surrogate boundary, this weight is increased by 20 in 2D and 40 in 3D for each
face that lies on $\tilde{\Gamma}$, accounting for the more expensive boundary term evaluations. The partitioning is
then performed by the parallel mesh distribution algorithms available in \texttt{deal.II}, which use these weights to
distribute cells among processors.

\subsection{CutFEM: A Comparative Framework}
\label{sec:cutfem_introduction}
For the same model Poisson problem $-\Delta u = f$ in $\Omega$ with $u = g$ on $\Gamma$, CutFEM seeks $u_h \in
    V_h^{\text{cut}}$ defined on the cut domain $\Omega_h = \bigcup_{K \cap \Omega \neq \emptyset} K \cap \Omega$,
where
the union is over all background mesh cells $K$ that intersect the physical domain $\Omega$. The CutFEM weak
formulation is:
\begin{align}
    \int_{\Omega_h} \nabla u_h \cdot \nabla v_h \, dx & - \int_{\Gamma_h} (\nabla u_h \cdot \mathbf{n}) v_h \, ds
    \nonumber
    \\
                                                      & - \int_{\Gamma_h} (\nabla v_h \cdot \mathbf{n}) (u_h - g) \, ds
    \nonumber
    \\
                                                      & + \int_{\Gamma_h} \sigma_{\Gamma} (u_h - g) v_h \, ds +
    s_{\text{GP}}(u_h, v_h) = \int_{\Omega_h} f v_h \, dx,
\end{align}
where $\Gamma_h$ is the discrete representation of the boundary $\Gamma$, and $s_{\text{GP}}(u_h, v_h)$ is the ghost
penalty stabilization term.

The ghost penalty~\cite{burman2010ghost,wichrowski2025matrix} is the key for CutFEM stability, as small cut elements
can lead to severe
ill-conditioning. It takes the form:
\begin{align}
    s_{\text{GP}}(u_h, v_h) = \sum_{F \in \mathcal{F}_{\text{GP}}}  \sum_{k=1,\ldots, p}	\gamma_F h_F
    \int_F \llbracket \frac{\partial^k
        u_h}{\partial^k n} \rrbracket
    \cdot \llbracket \frac{\partial^k v_h}{\partial^k n} \rrbracket \, ds,
\end{align}
where $\mathcal{F}_{\text{GP}}$ is the set of faces where the ghost penalty is applied (typically faces of cut cells
that are not on the boundary), $\llbracket \nabla u_h \rrbracket$ denotes the jump in the gradient across face $F$,
$h_F$ is the face diameter, and $\gamma_F > 0$ is a stabilization parameter.

The computational complexity of matrix-free CutFEM evaluation differs significantly from SBM due to irregular
integration domains and the ghost penalty mechanism. For cut cells, the complexity scales as $O(p^{2d})$ per cell due
to the inability to exploit sum factorization over arbitrarily shaped domains. The ghost penalty evaluation, however,
benefits from tensor-product structure~\cite{wichrowski2025matrix}: the gradient jumps $\llbracket \nabla u_h
    \rrbracket$
can be computed efficiently using a tensor product of 1D derivative matrices in one direction (across the face) and
$(d-1)$-dimensional mass matrices in the remaining directions. This yields a complexity of $O(p^{d+1})$ per face for
the ghost penalty. For our numerical comparisons, we utilize the matrix-free CutFEM implementation
from~\cite{wichrowski2025matrix}.

\section{Numerical Results}\label{sec:numerical_results}
This section presents numerical experiments demonstrating the computational efficiency of the matrix-free SBM operator
evaluation and comparing it with matrix-free CutFEM implementations. All computations are performed using a custom
implementation building upon the \texttt{deal.II} library~\cite{dealii2019design, kronbichler2012generic,
    bergbauer2024high}. The background mesh for all experiments consists of Cartesian cells.

The numerical experiments build upon the software framework developed in~\cite{wichrowski2025Geometric}. To validate
the correctness of the matrix-free operator evaluation, the resulting operators were compared against their sparse
matrix counterparts, ensuring exact reproducibility of matrix-vector products. The underlying finite element
operators are identical to those in~\cite{wichrowski2025Geometric}, we refer the reader to other works for detailed
convergence analysis of the SBM discretization~\cite{atallah2022high,collins2023penalty}.
In Appendix~\ref{app:convergence}, we provide a brief comparison of the convergence of SBM and CutFEM
discretizations for the Poisson problem. The focus here is exclusively on the efficient matrix-free evaluation of these
operators.

The experiments are structured to isolate different aspects of the matrix-free evaluation: local computational costs
through microbenchmarks, parallel scalability on realistic geometries, and the impact of geometric complexity on
operator performance. Comparisons with matrix-free CutFEM highlight the computational advantages of SBM's regular
integration domains.

All benchmarks were executed on a compute node equipped with two AMD EPYC 7282 16-core processors, providing a total of
32 physical cores. The implementation leverages AVX2 vectorization to maximize floating-point throughput, with a vector
width of 256 bits (four double-precision numbers).

For completeness, we note that traditional sparse matrix implementations of SBM operators exhibit significantly
inferior performance compared to the matrix-free approaches presented here. For reference, sparse matrix-vector
products for $p=3$ elements in 3D achieve throughput below $4 \times 10^7$ DoFs/sec for continuous Galerkin
discretizations, with even worse performance for discontinuous Galerkin formulations due to the increased coupling
between degrees of freedom. This represents more than an order of magnitude performance penalty compared to the
matrix-free implementations.

More critically, the matrix assembly phase required for sparse matrix approaches becomes prohibitively expensive for
high-order elements. The assembly process  must be performed whenever geometric or material properties change, making
it impractical for problems with evolving boundaries or nonlinear material behavior. The combination of expensive
assembly, large memory requirements (scaling as $O(p^{2d})$ per element), and poor matrix-vector product performance
makes sparse matrix approaches unsuitable for the high-performance computing applications that SBM is designed to
address. Consequently, we focus exclusively on matrix-free implementations throughout this work.

\subsection{Microbenchmarks: Local Computational Costs}\label{ssec:microbenchmarks}

The theoretical computational complexity of matrix-free SBM operator application scales as $O(N_c p^{d+1} + N_f p^d +
    N_{sf} p^{2d-1})$, where $N_c$ is the number of cells in $\tilde{\Omega}$, $N_f$ is the number of interior faces
(for
DG), $N_{sf}$ is the number of surrogate boundary faces, and $p$ is the polynomial degree. For large meshes, the
number of interior cells $N_c$ is much larger than the number of boundary faces $N_{sf}$, so the overall complexity is
dominated by the volumetric term, scaling as $O(N_c p^{d+1})$.
The most expensive component per face is the evaluation of shifted boundary terms, which scales as $O(p^{2d-1})$ due to
the inability to use sum factorization for evaluating functions at shifted points.

To validate this theoretical analysis and assess local computational costs, we perform microbenchmarks measuring the
wall time for applying operator contributions from individual mesh entities: interior cells, and surrogate boundary.
These measurements isolate the pure computational cost without memory bandwidth effects.
The benchmarks assume full utilization of vectorization. For the standard SBM cell, SBM face, and full matrix
evaluations, this is achieved by processing four identical entities simultaneously to leverage the CPU's AVX vector
width. In contrast, for a CutFEM cut cell, vectorization is performed over the quadrature points within that single
cell; its timing is therefore multiplied by four to provide a comparable throughput measurement.

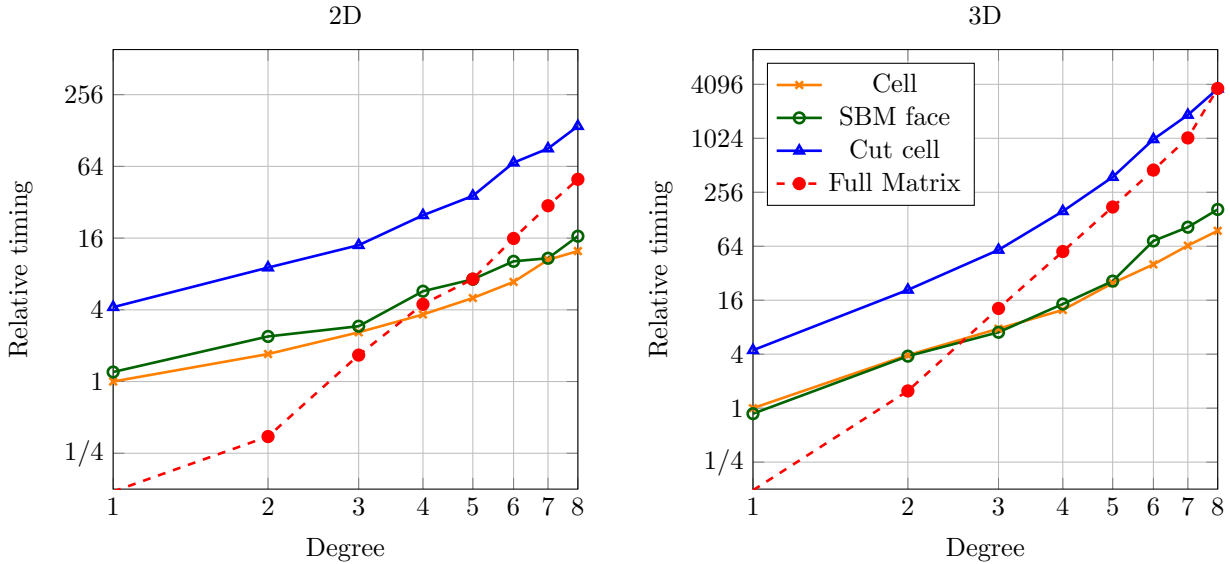
\begin{figure}[htbp]
    \pgfmathsetmacro{\twodnormal}{1.0480e-01} 
    \pgfmathsetmacro{\threednormal}{2.4230e-01} 

    \centering
    {
        \begin{tikzpicture}
            \begin{axis}[
                    title=2D,
                    xlabel=Degree,
                    ylabel=Relative timing,
                    legend pos=north west ,
                    width=0.47\textwidth,
                    height=0.45\textwidth,
                    grid=both,
                    ymode=log,
                    xmode=log,
                    log basis y=2,
                    ymin =0.125,
                    ymax=612,
                    xmin =1,
                    xmax = 8,
                    xtick={1,2,3,4,5,6,7,8},
                    xticklabels={1,2,3,4,5,6,7,8},
                    xticklabel style={/pgf/number format/fixed},
                    ytick={0.25, 1, 4, 16, 64, 256},
                    yticklabels={1/4, 1, 4, 16, 64, 256},
                ]

                \plotlinestyleone table [x index=0, y
                        index=2,  y expr=\thisrowno{2}/(\twodnormal), header=has
                        colnames]
                    {results/one_cell_2d.txt};

                \plotlinestyletwo table [x index=0, y
                        index=2,  y expr=\thisrowno{4}/(\twodnormal), header=has
                        colnames]
                    {results/one_cell_2d.txt};

                \plotlinestylethree table [x index=0, y
                        index=1, y expr=\thisrowno{1}*4/(\twodnormal), header=has
                        colnames]
                    {results/one_cell_2d.txt};

                \plotlinefull[dashed]
                table [x index=0, y
                        index=2,  y expr=\thisrowno{5}/(\twodnormal), header=has
                        colnames]
                    {results/one_cell_2d.txt};



            \end{axis}
        \end{tikzpicture}
    }
    \hfill
    {
        \begin{tikzpicture}
            \begin{axis}[
                    title=3D,
                    xlabel=Degree,
                    ylabel=Relative timing,
                    legend pos=north west,
                    width=0.47\textwidth,
                    height=0.45\textwidth,
                    grid=both,
                    ymode=log,
                    xmode=log,
                    log basis y=2,
                    ymin =0.125,
                    ymax = 0.1e5,
                    xmin =1,
                    xmax = 8,
                    xtick={1,2,3,4,5,6,7,8},
                    xticklabels={1,2,3,4,5,6,7,8},
                    ytick={0.25, 1, 4, 16, 64, 256, 1024, 4096, 16384},
                    yticklabels={1/4, 1, 4, 16, 64, 256,1024, 4096, 16384},
                ]
                \plotlinestyleone table [x index=0, y
                        index=2,  y expr=\thisrowno{2}/(\threednormal), header=has
                        colnames]
                    {results/one_cell_3d.txt};
                \addlegendentry{Cell}

                \plotlinestyletwo table [x index=0, y
                        index=2,  y expr=\thisrowno{4}/(\threednormal), header=has
                        colnames]
                    {results/one_cell_3d.txt};
                \addlegendentry{SBM face}

                \plotlinestylethree table [x index=0, y
                        index=1, y expr=\thisrowno{1}*4/(\threednormal), header=has
                        colnames]
                    {results/one_cell_3d.txt};
                \addlegendentry{Cut cell}

                \plotlinefull[dashed]
                table [x index=0, y
                        index=2,  y expr=\thisrowno{5}/(\threednormal), header=has
                        colnames]
                    {results/one_cell_3d.txt};
                \addlegendentry{Full Matrix}



            \end{axis}
        \end{tikzpicture}
    }
    \caption{Microbenchmark:  Relative time per application for different evaluation methods on a single cell.
        The time is	normalized by the time of FEEvaluation for $k=1$, that is 0.105 $\mu$s in 2D and 0.2423 $\mu$s
        in 3D.
        Expected single element/face evaluation times for SBM (CG/DG) and CutFEM for varying
        polynomial degree $k_0$ in $d=2$ or $d=3$ dimensions. SBM evaluations on regular grid entities are anticipated
        to be
        faster than CutFEM evaluations on arbitrarily cut cells, especially as $k_0$ increases, due to simpler
        integration
    }\label{fig:microbenchmark}
\end{figure}

Figure~\ref{fig:microbenchmark} presents timing results comparing SBM and CutFEM local evaluations across polynomial
degrees $p=1$ to $p=8$ in both 2D and 3D. The results demonstrate that SBM local evaluations are consistently faster
than their CutFEM counterparts, with the performance advantage becoming more pronounced at higher polynomial degrees.
It is important to note that the reported CutFEM timings are based on a standard Gauss quadrature rule applied to the
cut cell geometry. In practice, more sophisticated and computationally expensive quadrature rules are often necessary
to accurately integrate over the irregular domains of cut cells, which would further increase the computational
workload for CutFEM compared to SBM.

It is also instructive to compare these matrix-free techniques with a \emph{full matrix} evaluation, where
the local cell matrix is pre-assembled and applied via a matrix-vector product. While this approach has the highest
asymptotic complexity of $O(N_c^2)$, for lower polynomial degrees it can be faster than sum factorization, as seen in
the benchmarks. This is due to highly optimized linear algebra libraries and the assumed vectorization across
multiple identical cells, a strategy that is possible in hybrid matrix-free methods where only geometry-affected cells
are assembled. Notably, this assembled approach still outperforms the evaluation on a CutFEM cut cell. However, this
performance comes at the significant upfront cost of matrix assembly and storage. The assembly itself requires $N_c$
cell-operator-vector products, each with a cost of $O(p^{2d})$, leading to a total complexity of $O(N_c p^{2d})$. This
assembly cost can be prohibitive, especially for high polynomial degrees $p$. Furthermore, the large memory footprint
and the memory-bandwidth-bound nature of the subsequent matrix-vector products negate the primary advantages of
matrix-free methods for high-order or large-scale computations. For vector-valued problems with multiple components
(such as elasticity or fluid dynamics), the situation becomes even worse, as matrix sizes grow quadratically with the
number of components, making full matrix approaches particularly inefficient for such multi-physics applications.

\begin{figure}[htbp]

    \centering
    {
        \begin{tikzpicture}
            \begin{axis}[
                    title=2D,
                    xlabel=Degree,
                    ylabel=Memory usage,
                    legend pos=north west ,
                    width=0.47\textwidth,
                    height=0.45\textwidth,
                    grid=both,
                    ymode=log,
                    xmode=log,
                    log basis y=2,
                    xmin =1,
                    xmax = 8,
                    xtick={1,2,3,4,5,6,7,8},
                    xticklabels={1,2,3,4,5,6,7,8},
                    xticklabel style={/pgf/number format/fixed},
                    ytick={0.25, 1, 4, 16, 64, 256, 1024, 4096, 16384},
                    yticklabels={1/4, 1, 4, 16, 64, 256,1024, 4096, 16384},
                ]

                \plotlinestyletwo table [x index=0, y
                        index=2,  y expr=(\thisrowno{2}- \thisrowno{4})/8, header=has
                        colnames]
                    {results/one_cell_mem_2d.txt};

                \plotlinestylethree table [x index=0, y
                        index=1, y expr=(\thisrowno{1}- \thisrowno{4})/8, header=has
                        colnames]
                    {results/one_cell_mem_2d.txt};

                \plotlinefull[dashed] table [x index=0, y
                        index=3,  y expr=(\thisrowno{3})/8/4, header=has
                        colnames]
                    {results/one_cell_mem_2d.txt};



            \end{axis}
        \end{tikzpicture}
    }
    \hfill
    {
        \begin{tikzpicture}
            \begin{axis}[
                    title=3D,
                    xlabel=Degree,
                    ylabel=Memory usage,
                    legend pos=north west,
                    width=0.47\textwidth,
                    height=0.45\textwidth,
                    grid=both,
                    ymode=log,
                    xmode=log,
                    log basis y=2,
                    xmin =1,
                    xmax = 8,
                    xtick={1,2,3,4,5,6,7,8},
                    xticklabels={1,2,3,4,5,6,7,8},
                    ytick={0.25, 1, 4, 16, 64, 256, 1024, 4096, 16384, 65536, 262144},
                    yticklabels={1/4, 1, 4, 16, 64, 256,1024, 4096, 16k, 65k, 262k},
                ]

                \plotlinestyletwo table [x index=0, y
                        index=2,  y expr=(\thisrowno{2}- \thisrowno{4})/8, header=has
                        colnames]
                    {results/one_cell_mem_3d.txt};
                \addlegendentry{SBM face}

                \plotlinestylethree table [x index=0, y
                        index=1, y expr=(\thisrowno{1}- \thisrowno{4})/8, header=has
                        colnames]
                    {results/one_cell_mem_3d.txt};
                \addlegendentry{Cut cell}

                \plotlinefull[dashed] table [x index=0, y
                        index=3,  y expr=(\thisrowno{3})/8/4, header=has
                        colnames]
                    {results/one_cell_mem_3d.txt};
                \addlegendentry{Full Matrix}



            \end{axis}
        \end{tikzpicture}
    }
    \caption{Microbenchmark: Memory consumption for different evaluation methods on a single cell/face.
        The memory usage is reported as the number of stored in double precision values.
        The plot compares the memory required for precomputed data for an SBM face and a CutFEM cut cell against
        the storage for a full local matrix, for varying polynomial degree $k$ in 2D and 3D.
        The memory for SBM and CutFEM is only for the additional data structures needed for matrix-free evaluation on
        unfitted geometries.
    }\label{fig:microbenchmark_memory}
\end{figure}
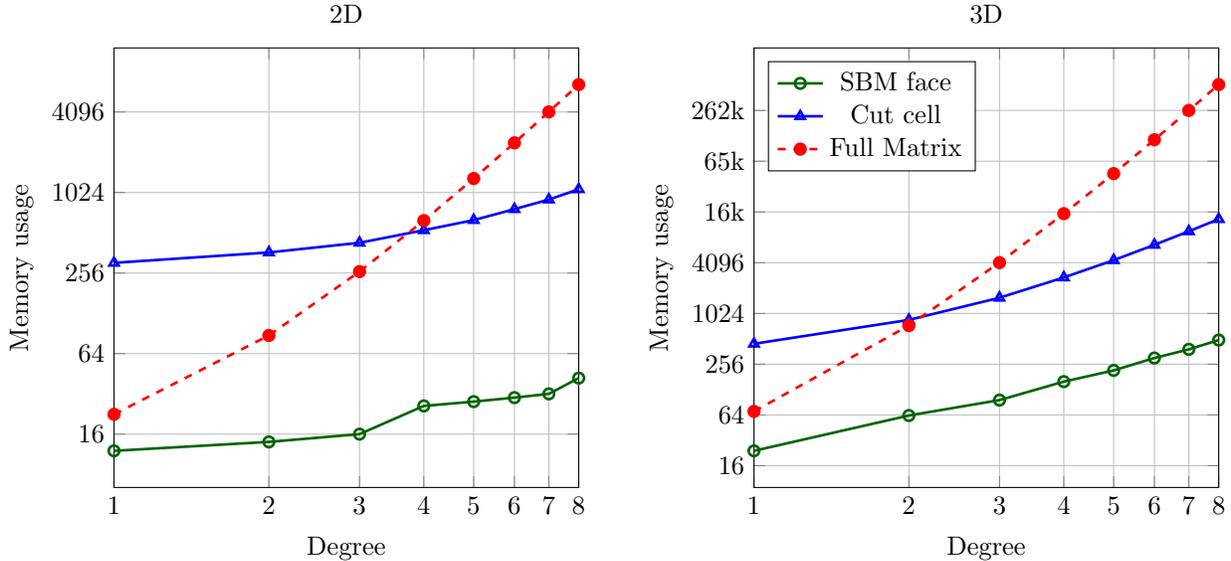

Figure~\ref{fig:microbenchmark_memory} compares memory requirements for precomputed geometric data.
We do not include the memory for the interior cell data, as this can be reused across all cells in the background mesh
and is therefore independent of the number of cells. All computations are performed in double precision, and memory
usage is
reported in units of double precision numbers (bytes divided by 8). SBM requires significantly less memory per
quadrature point on the surrogate boundary, avoiding the complex cut cell information and adaptive quadrature data
structures needed by CutFEM.
This memory efficiency translates to better cache utilization and reduced data movement costs during operator
evaluation.

\subsection{Parallel Performance: Single Ball Benchmark}\label{ssec:performance_analysis}

To evaluate the parallel scalability and real-world performance of the matrix-free SBM operator, we solve a Poisson
problem on a unit ball domain using both Continuous Galerkin (CG) and Discontinuous Galerkin (DG) discretizations. This
benchmark uses 32 MPI ranks on a single compute node, representing a typical high-performance computing setup where
memory bandwidth and inter-process communication are optimized.

The test measures the total time for a fixed number of matrix-vector products (representative of an iterative solver's
computational kernel) and computes throughput in terms of degrees of freedom processed per second. This metric captures
the essential performance characteristic for iterative solvers, where the matrix-vector product typically dominates the
computational cost.

\begin{figure}[htbp]
    \pgfmathsetmacro{\ncycles}{10} 
    \centering
    {
        \begin{tikzpicture}
            \begin{axis}[
                    title=2D,
                    xlabel={Number of active cells},
                    xlabel style={yshift=-1em},
                    ylabel={Throughput (DoFs/sec)},
                    legend pos=south east,
                    width=0.47\textwidth,
                    height=0.45\textwidth,
                    grid=both,
                    ymode=log,
                    xmode=log,
                    log basis y=10,
                    log basis x=10,
                    xmin=50,
                    xmax=4663636,
                    xtick={52, 256,1060,4416,17940,72340,290508,1164688,4663636},
                    xticklabels={52, 256,1.1k,4.4k,18k,72k,291k,1.16M,4.66M},
                    xticklabel style={rotate=45,anchor=east},
                ]
                \plotlinestyleone table [
                        x expr=\thisrowno{2},
                        y expr={ \thisrowno{8}},
                        col sep=space,
                        skip first n=1
                    ] {results/single_ball_CG_p1_d2.txt};

                \plotlinestyletwo table [
                        x expr=\thisrowno{2},
                        y expr={ \thisrowno{8}},
                        col sep=space,
                        skip first n=1
                    ] {results/single_ball_CG_p2_d2.txt};

                \plotlinestylethree table [
                        x expr=\thisrowno{2},
                        y expr={ \thisrowno{8}},
                        col sep=space,
                        skip first n=1
                    ] {results/single_ball_CG_p3_d2.txt};

                \plotlinestyleone[dashed] table [
                        x expr=\thisrowno{2},
                        y expr={\thisrowno{2}*4 /(\thisrowno{7}*1e-6 / \ncycles)},
                        col sep=space,
                        skip first n=1
                    ] {results/single_ball_p1_d2.txt};

                \plotlinestyletwo[dashed] table [
                        x expr=\thisrowno{2},
                        y expr={\thisrowno{2}*9 /(\thisrowno{7}*1e-6/ \ncycles)},
                        col sep=space,
                        skip first n=1
                    ] {results/single_ball_p2_d2.txt};

                \plotlinestylethree[dashed] table [
                        x expr=\thisrowno{2},
                        y expr={ \thisrowno{8}},
                        col sep=space,
                        skip first n=1
                    ] {results/single_ball_p3_d2.txt};

            \end{axis}
        \end{tikzpicture}  }
    \hfill
    {
        \begin{tikzpicture}
            \begin{axis}[
                    title=3D,
                    xlabel={Number of active cells},
                    xlabel style={yshift=-1em},
                    ylabel={Throughput (DoFs/sec)},
                    legend pos=south east,
                    width=0.47\textwidth,
                    height=0.45\textwidth,
                    grid=both,
                    ymode=log,
                    xmode=log,
                    log basis y=10,
                    log basis x=10,
                    ymin=4e5,
                    ymax=2e9,
                    xmin=8,
                    xmax=913936,
                    xtick={8, 136, 1352, 12688, 109696, 913936},
                    xticklabels={8, 136, 1.4k, 12.7k, 110k, 914k},
                    xticklabel style={rotate=45,anchor=east},	   ]

                \plotlinestyleone table [
                        x expr=\thisrowno{2},
                        y expr={ \thisrowno{9}},
                        col sep=space,
                        skip first n=1
                    ] {results/single_ball_CG_p1_d3.txt};
                \addlegendentry{$p=1$}

                \plotlinestyletwo table [
                        x expr=\thisrowno{2},
                        y expr={ \thisrowno{9}},
                        col sep=space,
                        skip first n=1
                    ] {results/single_ball_CG_p2_d3.txt};
                \addlegendentry{$p=2$}

                \plotlinestylethree table [
                        x expr=\thisrowno{2},
                        y expr={ \thisrowno{9}},
                        col sep=space,
                        skip first n=1
                    ] {results/single_ball_CG_p3_d3.txt};
                \addlegendentry{$p=3$}

                \plotlinestyleone[dashed] table [
                        x expr=\thisrowno{2},
                        y expr={ \thisrowno{9}},
                        col sep=space,
                        skip first n=1
                    ] {results/single_ball_p1_d3.txt};

                \plotlinestyletwo[dashed] table [
                        x expr=\thisrowno{2},
                        y expr={ \thisrowno{9}},
                        col sep=space,
                        skip first n=1
                    ] {results/single_ball_p2_d3.txt};

                \plotlinestylethree[dashed] table [
                        x expr=\thisrowno{2},
                        y expr={ \thisrowno{9}},
                        col sep=space,
                        skip first n=1
                    ] {results/single_ball_p3_d3.txt};

            \end{axis}
        \end{tikzpicture}  }
    \caption{Throughput in degrees of freedom per second (DoFs/sec) for the SBM operator on a unit ball geometry, shown
        for 2D (left) and 3D (right) computations. The plots compare the performance for different polynomial degrees
        $p=1, 2,  3$ as a function of the number of active cells. Solid lines represent results for continuous
        elements, while discontinuous lines correspond to discontinuous Galerkin (DG) elements.
    }\label{fig:unit_ball_performance}
\end{figure}
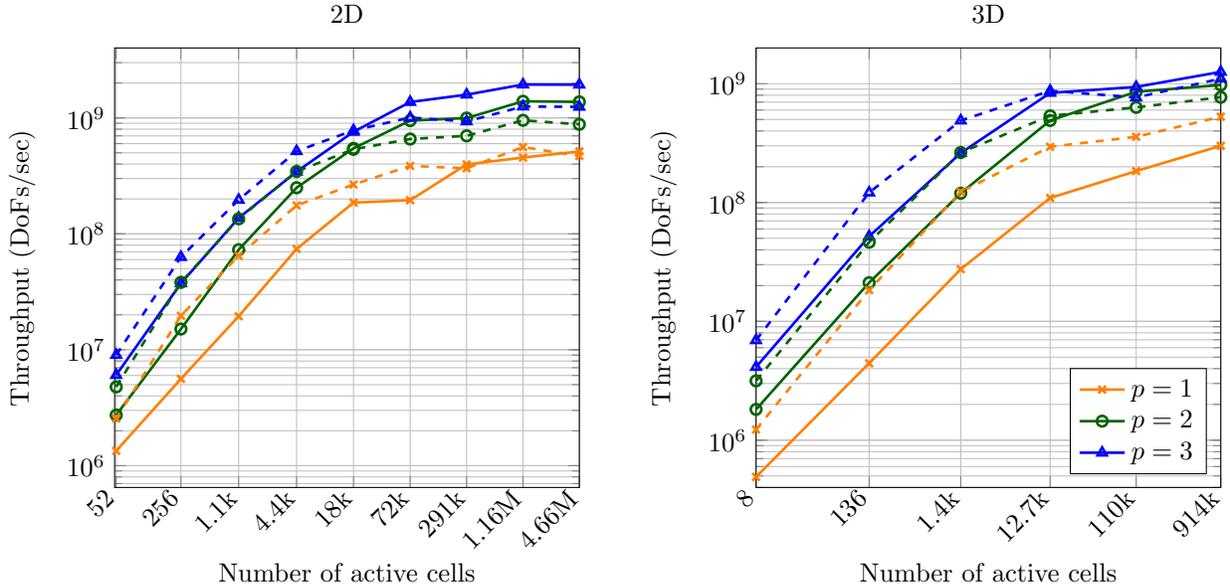

Figure~\ref{fig:unit_ball_performance} demonstrates the scaling characteristics of the method across polynomial degrees
$p=1$ to $p=3$. The results showcase high throughput, particularly for higher polynomial degrees, highlighting the
benefits of matrix-free techniques where the computational cost per degree of freedom can decrease with increasing $p$
due to improved arithmetic intensity and better cache utilization from sum factorization.

For DG discretizations, the additional interior face terms contribute to the computational cost but maintain good
scalability due to the effective use of sum factorization on the regular face quadrature points. The performance
characteristics demonstrate that both CG and DG variants of SBM benefit significantly from the matrix-free approach,
especially at higher polynomial degrees.

\subsection{Impact of Geometric Complexity: Multiple Ball Benchmark}\label{ssec:multiple_balls}

To investigate how geometric complexity affects operator evaluation performance, we conduct tests using domains
containing multiple randomly placed balls within the background mesh. The number of balls is varied, starting from
three
large balls that create a domain with a relatively few intersections, and increasing to 25 smaller
balls. In the 25-ball case, the geometry becomes highly complex, with the fraction of intersected cells reaching
approximately 50\%. This setup allows for a systematic study of performance as the
boundary-to-volume ratio of the unfitted geometry increases.

These tests employ a specialized execution strategy designed to isolate pure operator evaluation costs from parallel
overhead effects. While running in serial mode, 32 matrix-vector products execute simultaneously in perfect
synchronization, fully utilizing CPU resources and memory bandwidth.
By using the emulated parallel
approach in the multiple ball tests, we bypass these practical constraints, allowing us to observe the pure
computational performance potential of the kernel without the limitations imposed by distributed computing
environments.  To make the results comparable
with the fully parallel benchmarks in Figure~\ref{fig:unit_ball_performance}, the resulting single-core throughput is
multiplied by the number of cores (32) to estimate the total achievable throughput.

The background mesh consists of $64^3 = 262144$ cells,	 with the number of intersected cells varying as the
geometric complexity increased.  As the number of balls increases, the surrogate domain shrinks, with the number of
active cells decreasing to 15,078, while the number of intersected cells peaks at 23,363 for 10 balls.	For both SBM
and CutFEM, the fraction of intersected cells is defined as the number of intersected cells divided by the total number
of cells that are either active or intersected, ensuring a consistent geometric complexity metric across both methods.
In case of 25 balls, the fraction of intersected cells reaches approximately 50\%.

\begin{figure}[htbp]

    \pgfmathsetmacro{\ncores}{32} 
    \pgfmathsetmacro{\nskip}{3} 

    \centering
    {
        \begin{tikzpicture}
            \begin{axis}[
                    xlabel={Fraction of intersected cells},
                    ylabel={Throughput (DoFs/sec)},
                    legend pos=south west,
                    width=0.47\textwidth,
                    height=0.45\textwidth,
                    grid=both,
                    ymode=log,
                    log basis y=10,
                    xticklabel style={rotate=45,anchor=east},
                    xmin=0,
                    xmax=0.51,
                ]
                \plotlinestyleone
                table [
                        x expr=\thisrowno{1}/(\thisrowno{0}+\thisrowno{1}),
                        y expr={ \thisrowno{10}*\ncores},
                        col sep=space,
                        skip first n=\nskip]
                    {results/fake_parallel/3D_p1.txt};

                \plotlinestyletwo
                table [
                        x expr=\thisrowno{1}/(\thisrowno{0}+\thisrowno{1}),
                        y expr={ \thisrowno{10}*\ncores},
                        col sep=space,
                        skip first n=\nskip
                    ] {results/fake_parallel/3D_p2.txt};

                \plotlinestylethree
                table [
                        x expr=\thisrowno{1}/(\thisrowno{0}+\thisrowno{1}),
                        y expr={ \thisrowno{10}*\ncores},
                        col sep=space,
                        skip first n=\nskip
                    ] {results/fake_parallel/3D_p3.txt};

                \plotlinestyleone[dashed]
                table [
                        x expr=\thisrowno{11},
                        y expr={\thisrowno{3} *\ncores / \thisrowno{4}},
                        col sep=space,
                        skip first n=\nskip
                    ] {results/cutFEMref/fake_parallel_balls/3D_p1.txt};

                \plotlinestyletwo[dashed]
                table [
                        x expr=\thisrowno{11},
                        y expr={\thisrowno{3} *\ncores / \thisrowno{4}},
                        col sep=space,
                        skip first n=\nskip
                    ] {results/cutFEMref/fake_parallel_balls/3D_p2.txt};

                \plotlinestylethree[dashed]
                table [
                        x expr=\thisrowno{11},
                        y expr={\thisrowno{3} *\ncores / \thisrowno{4}},
                        col sep=space,
                        skip first n=\nskip
                    ] {results/cutFEMref/fake_parallel_balls/3D_p3.txt};

            \end{axis}
        \end{tikzpicture}}
    \hfill
    {
        \begin{tikzpicture}
            \begin{axis}[
                    xlabel={Fraction of intersected cells},
                    ylabel={Geometry time fraction},
                    legend pos=south east,
                    width=0.47\textwidth,
                    height=0.45\textwidth,
                    grid=both,
                    xticklabel style={rotate=45,anchor=east},
                    xmin=0,
                    xmax=0.51,
                    ytick={ 0.25, 0.5, 0.75, 1},
                    yticklabels={ 25\%, 50\%, 75\%, 100\%},
                    ymax=1,
                    ymin=0,
                ]

                \plotlinestyleone
                table [
                        x expr=\thisrowno{1}/(\thisrowno{0}+\thisrowno{1}),
                        y expr={ \thisrowno{13}/ \thisrowno{9}},
                        col sep=space,
                        skip first n=\nskip
                    ] {results/fake_parallel/3D_p1.txt};
                \addlegendentry{$p=1$}

                \plotlinestyletwo
                table [
                        x expr=\thisrowno{1}/(\thisrowno{0}+\thisrowno{1}),
                        y expr={ \thisrowno{13}/ \thisrowno{9}},
                        col sep=space,
                        skip first n=\nskip
                    ] {results/fake_parallel/3D_p2.txt};
                \addlegendentry{$p=2$}

                \plotlinestylethree
                table [
                        x expr=\thisrowno{1}/(\thisrowno{0}+\thisrowno{1}),
                        y expr={ \thisrowno{13}/ \thisrowno{9}},
                        col sep=space,
                        skip first n=\nskip
                    ] {results/fake_parallel/3D_p3.txt};
                \addlegendentry{$p=3$}

                \plotlinestyleone[dashed]
                table [
                        x expr=\thisrowno{11},
                        y expr={(\thisrowno{6}+\thisrowno{7}) / \thisrowno{4}},
                        col sep=space,
                        skip first n=\nskip
                    ] {results/cutFEMref/fake_parallel_balls/3D_p1.txt};

                \plotlinestyletwo[dashed]
                table [
                        x expr=\thisrowno{11},
                        y expr={(\thisrowno{6}+\thisrowno{7}) / \thisrowno{4}},
                        col sep=space,
                        skip first n=\nskip
                    ] {results/cutFEMref/fake_parallel_balls/3D_p2.txt};

                \plotlinestylethree[dashed]
                table [
                        x expr=\thisrowno{11},
                        y expr={(\thisrowno{6}+\thisrowno{7}) / \thisrowno{4}},
                        col sep=space,
                        skip first n=\nskip
                    ] {results/cutFEMref/fake_parallel_balls/3D_p3.txt};

            \end{axis}
        \end{tikzpicture}}
    \caption{Performance analysis for multiple ball benchmark as a function of fraction of intersected cells.
        Left: throughput in degrees of freedom processed per second.
        Right: fraction of operator evaluation time spent on geometry-related computations.
        Solid lines indicate SBM results, dashed lines indicate CutFEM results.
    }\label{fig:multiple_ball_performance}
\end{figure}
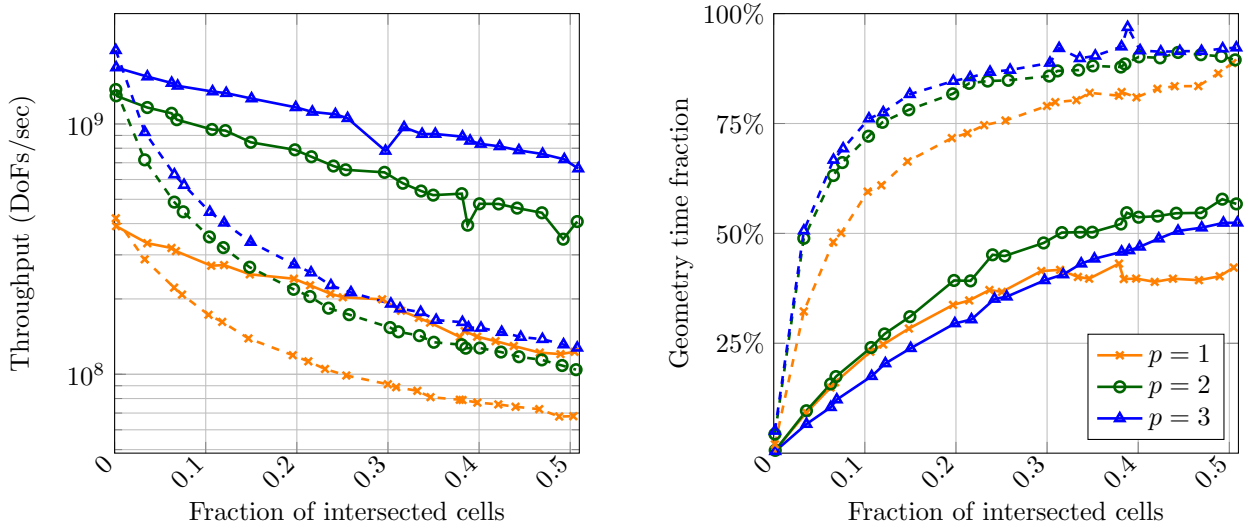

Figure~\ref{fig:multiple_ball_performance} presents both throughput measurements and the fraction of time spent on
geometry-related computations as the number of balls (and correspondingly, the fraction of intersected cells)
increases. The left panel shows that SBM maintains efficient operator evaluation even as geometric complexity grows,
with only moderate throughput reduction for high fractions of intersected cells.  The results also indicate that even
with a high fraction of
intersected cells, the SBM throughput remains competitive, dropping by less than one order of
magnitude from the ideal case.

It is worth noting that the maximum throughput achieved in this benchmark for $p=3$ with minimal
intersected cells is
approximately $1.68\times 10^9$
DoFs/sec, which is about 26\% higher than the $1.33\times 10^9$ DoFs/sec achieved in the unit ball benchmark under full
MPI parallelization. This difference highlights the performance costs of communication overhead and load balancing
challenges present in the fully parallel unit ball benchmark.

The right panel quantifies the computational overhead specifically attributable to geometry-related operations,
revealing a stark contrast between the two methods. For SBM, geometry-related operations refer exclusively to the
evaluation of surrogate boundary faces, while for CutFEM, it encompasses both the processing of intersected cells and
the evaluation of ghost penalty terms. The results show that evaluation of CutFEM operator quickly becomes saturated by
geometry-related computations, with this fraction rapidly approaching 80--90\% as the domain complexity increases. In
contrast, even with half of cells being intersected, the geometry-related computational overhead for SBM remains
manageable, staying below 60\% across all polynomial degrees. This confirms SBM's resilience to geometric
complexity, as first predicted by the microbenchmark results in Figure~\ref{fig:microbenchmark}.
The regular integration domains in the surrogate boundary approach permit more efficient evaluation
patterns compared to the irregular cut cells and stabilization terms in CutFEM.

\subsection{Initialization of Operator and Memory Requirements}\label{ssec:initialization_performance}

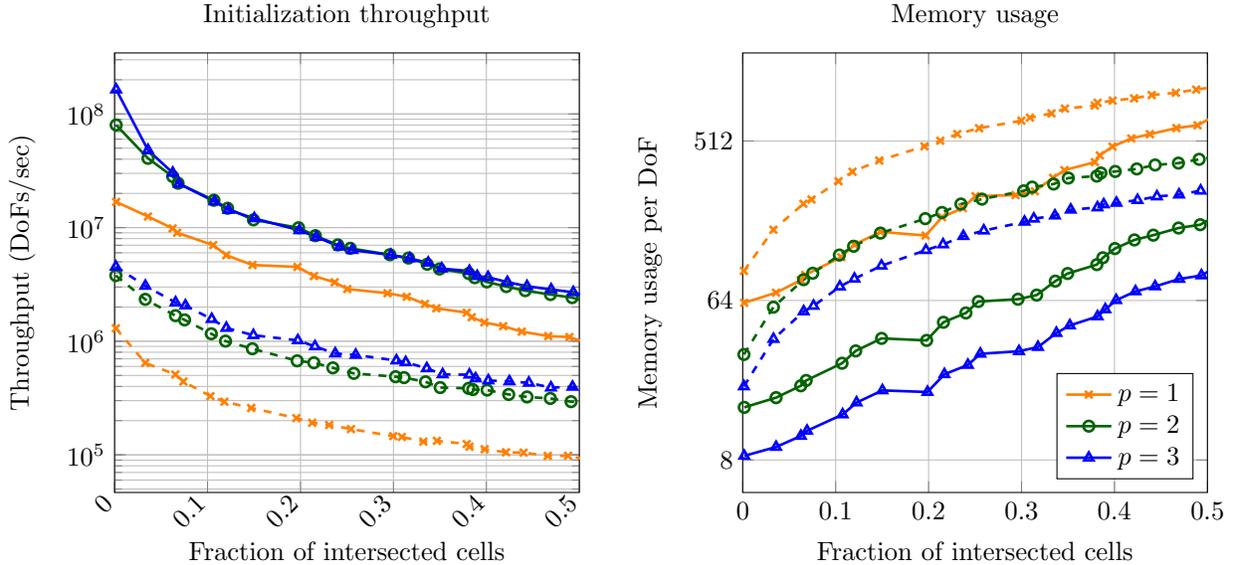
\begin{figure}[htbp]

    \pgfmathsetmacro{\ncores}{32} 
    \pgfmathsetmacro{\nskip}{3} 

    \centering
    {
        \begin{tikzpicture}
            \begin{axis}[
                    title={Initialization throughput},
                    xlabel={Fraction of intersected cells},
                    ylabel={Throughput (DoFs/sec)},
                    legend pos=north east,
                    width=0.47\textwidth,
                    height=0.45\textwidth,
                    grid=both,
                    ymode=log,
                    log basis y=10,
                    xticklabel style={rotate=45,anchor=east},
                    xmin=0,
                    xmax=0.5,
                ]

                \plotlinestyleone
                table [
                        x expr=\thisrowno{1}/(\thisrowno{0}+\thisrowno{1}),
                        y expr={  \thisrowno{4}  \ncores * 1e6	/ \thisrowno{6} },
                        col sep=space,
                        skip first n=\nskip
                    ] {results/fake_parallel/3D_p1.txt};

                \plotlinestyletwo
                table [
                        x expr=\thisrowno{1}/(\thisrowno{0}+\thisrowno{1}),
                        y expr={  \thisrowno{4}  \ncores * 1e6	/ \thisrowno{6} },
                        col sep=space,
                        skip first n=\nskip
                    ] {results/fake_parallel/3D_p2.txt};

                \plotlinestylethree
                table [
                        x expr=\thisrowno{1}/(\thisrowno{0}+\thisrowno{1}),
                        y expr={  \thisrowno{4}  \ncores * 1e6	/ \thisrowno{6} },
                        col sep=space,
                        skip first n=\nskip
                    ] {results/fake_parallel/3D_p3.txt};

                \plotlinestyleone[dashed]
                table [
                        x expr=\thisrowno{11},
                        y expr={ \thisrowno{3}/ \thisrowno{10} },
                        col sep=space,
                        skip first n=\nskip
                    ] {results/cutFEMref/fake_parallel_balls/3D_p1.txt};

                \plotlinestyletwo[dashed]
                table [
                        x expr=\thisrowno{11},
                        y expr={ \thisrowno{3}/ \thisrowno{10} },
                        col sep=space,
                        skip first n=\nskip
                    ] {results/cutFEMref/fake_parallel_balls/3D_p2.txt};

                \plotlinestylethree[dashed]
                table [
                        x expr=\thisrowno{11},
                        y expr={ \thisrowno{3}/ \thisrowno{10} },
                        col sep=space,
                        skip first n=\nskip
                    ] {results/cutFEMref/fake_parallel_balls/3D_p3.txt};

            \end{axis}
        \end{tikzpicture}}
    \hfill
    {
        \begin{tikzpicture}
            \begin{axis}[
                    title={Memory usage},
                    xlabel={Fraction of intersected cells},
                    ylabel={Memory usage  per DoF},
                    legend pos=south east,
                    width=0.47\textwidth,
                    height=0.45\textwidth,
                    grid=both,
                    ymode=log,
                    log basis y=2,
                    ytick={8, 64, 512, 4096},
                    yticklabels={8, 64, 512, 4096},
                    yminorgrids=true,
                    ymajorgrids=true,
                    xmin=0,
                    xmax=0.5,
                ]

                \plotlinestyleone
                table [
                        x expr=\thisrowno{1}/(\thisrowno{0}+\thisrowno{1}),
                        y expr={ \thisrowno{16} /8 / \thisrowno{4} },
                        col sep=space,
                        skip first n=\nskip
                    ] {results/fake_parallel/3D_p1.txt};
                \addlegendentry{$p=1$}

                \plotlinestyletwo
                table [
                        x expr=\thisrowno{1}/(\thisrowno{0}+\thisrowno{1}),
                        y expr={ \thisrowno{16} / 8/ \thisrowno{4} },
                        col sep=space,
                        skip first n=\nskip
                    ] {results/fake_parallel/3D_p2.txt};
                \addlegendentry{$p=2$}

                \plotlinestylethree
                table [
                        x expr=\thisrowno{1}/(\thisrowno{0}+\thisrowno{1}),
                        y expr={ \thisrowno{16}  /8 / \thisrowno{4} },
                        col sep=space,
                        skip first n=\nskip
                    ] {results/fake_parallel/3D_p3.txt};
                \addlegendentry{$p=3$}

                \plotlinestyleone[dashed]
                table [
                        x expr=\thisrowno{11},
                        y expr={ \thisrowno{12} /\thisrowno{3}},
                        col sep=space,
                        skip first n=\nskip
                    ] {results/cutFEMref/fake_parallel_balls/3D_p1.txt};

                \plotlinestyletwo[dashed]
                table [
                        x expr=\thisrowno{11},
                        y expr={ \thisrowno{12} /\thisrowno{3}},
                        col sep=space,
                        skip first n=\nskip
                    ] {results/cutFEMref/fake_parallel_balls/3D_p2.txt};

                \plotlinestylethree[dashed]
                table [
                        x expr=\thisrowno{11},
                        y expr={ \thisrowno{12} /\thisrowno{3}},
                        col sep=space,
                        skip first n=\nskip
                    ] {results/cutFEMref/fake_parallel_balls/3D_p3.txt};

            \end{axis}
        \end{tikzpicture}  }
    \caption{Performance comparison between SBM (solid lines) and CutFEM (dashed lines) as a function of geometric
        complexity (fraction of intersected cells). Left: Initialization throughput in degrees of freedom (DoFs) per
        second.
        Right: Memory requirements per degree of freedom (in units of double precision numbers).
    }\label{fig:initializing}
\end{figure}

While operator evaluation performance is crucial, the initialization phase of unfitted methods significantly impacts
their practicality, especially for problems with evolving geometries. The left panel of Figure~\ref{fig:initializing}
compares the initialization throughput of SBM and CutFEM methods. SBM initialization involves identifying active cells
and computing closest point projections from surrogate boundary faces to the true boundary. Using a level set
representation, this projection is a straightforward, parallelizable Newton solve. In contrast, CutFEM initialization
is more complex, requiring the computation of cell-domain intersections, the generation of adaptive quadrature rules
for each cut cell, and the construction of ghost penalty connectivity; these geometric operations are inherently more
irregular and computationally intensive.

To maintain a consistent performance metric across all computational stages, we measure initialization throughput in
degrees of freedom per second (DoFs/sec). While the initialization process operates on mesh entities --- cells for
CutFEM
and faces for SBM --- its purpose is to construct the data structures required for all degrees of freedom. This metric
enables a direct comparison between initialization and operator evaluation throughput
(Figures~\ref{fig:unit_ball_performance} and~\ref{fig:multiple_ball_performance}), providing a holistic view of each
method's performance.

SBM demonstrates a substantially higher initialization throughput --- often by an order of magnitude --- compared to
CutFEM
across all tested scenarios. For typical problem sizes, the SBM initialization time is approximately equivalent to 10
matrix-vector products, making it a relatively modest overhead compared to the tens of iterations typically required in
well-preconditioned iterative solvers (though more iterations may be needed due to the non-matching nature of unfitted
methods). While the initialization throughput for both methods predictably declines as geometric complexity increases,
they do so at a comparable rate.

The right panel of Figure~\ref{fig:initializing} compares the memory requirements per degree of freedom for the
geometric data structures of each method, confirming the predictions from the microbenchmarks in
Figure~\ref{fig:microbenchmark_memory}. All computations are performed in double precision, and memory usage is
reported in units of double precision numbers (bytes divided by 8). SBM's memory advantage stems from its simple data
storage, requiring only the
shift vector and reference coordinates for each boundary quadrature point. In contrast, CutFEM requires larger data
structures for each cut cell, including custom quadrature. Consequently, CutFEM's memory usage per DoF increases with
geometric complexity as more cells are cut.

This memory advantage becomes more pronounced for larger problems, as the amount of boundary data scales more favorably
than the volumetric data. The lower memory footprint of SBM is crucial for modern hardware, where performance is often
limited not just by memory capacity but also by memory bandwidth bottlenecks.

\section{Conclusion and Future Work}\label{sec:conclusion}

We have presented a comprehensive matrix-free framework for the Shifted Boundary Method, applicable to both Continuous
Galerkin (CG) and Discontinuous Galerkin finite element discretizations. By leveraging tensor-product structures
for sum factorization on interior cells and faces, and employing efficient evaluation techniques for surrogate boundary
terms, the method achieves high computational performance and scalability, particularly for high-order elements.

The numerical experiments demonstrate significant computational advantages over matrix-free CutFEM implementations,
including faster local operator evaluations, reduced memory requirements, and better scalability with increasing
geometric complexity. These benefits stem from SBM's regular integration domains, in contrast to CutFEM's need for
complex quadrature rules and irregular data structures for arbitrarily cut cells.

The matrix-free operator evaluation provides an essential computational foundation for scalable iterative solvers
applied to large-scale unfitted finite element problems. While effective preconditioning strategies remain an
open challenge, the efficient operator application presented here is the key enabler for future developments in
efficient SBM solution strategies.

\section*{Acknowledgements}
The author would like to thank Guido Kanschat and Guglielmo Scovazzi for insightful discussions. The author is also
grateful to Luca Heltai for the suggestion to compare SBM against CutFEM; the author hopes the resulting
analysis resolves the discussion.

The author declares the use of language models (ChatGPT, Gemini, and Claude) to improve the clarity and readability of
the manuscript. All scientific content and technical claims are solely the responsibility of the author.

\bibliographystyle{siam}
\bibliography{literature} 

\newpage
\appendix
\section{Convergence Comparison: SBM vs. CutFEM}\label{app:convergence}
For complete comparison of SBM and CutFEM, we analyze the accuracy of the SBM and CutFEM implementations. We solve the
Poisson problem $-\Delta u = f$
on a domain $\Omega$ defined as a unit circle centered at the origin. The manufactured solution is given by:
\[
    u(\mathbf{x}) = 2 \cos(x_1) \sin(x_2).
\]
The right-hand side $f = -\Delta u$ and the Dirichlet boundary data $g = u|_{\partial\Omega}$ are derived from this
exact solution. We use a sequence of uniformly refined background meshes and compute the numerical solution for
polynomial degrees $p=1, 2$, and $3$. The error is measured in the $L^2$ norm over the respective computational domain
($\tilde{\Omega}$ for SBM and the cut domain $\Omega_h$ for CutFEM).

\begin{figure}[!ht]
    \centering
    \begin{tikzpicture}
        \begin{axis}[
                xlabel={$h$},
                ylabel={$L^2$ Error},
                xmode=log,
                ymode=log,
                log basis x=10,
                log basis y=10,
                legend pos=south east,
                width=0.47\textwidth,
                height=0.45\textwidth,
                grid=both,
            ]

            \plotlinestyleone
            table [x index=1, y index=2, col sep=space, skip first n=1]
                {results/convergence/sbm_p1.txt};
            \addlegendentry{$p=1$}
            \plotlinestyletwo table[x index=1, y index=2, col sep=space, skip first n=1]
                {results/convergence/sbm_p2.txt};
            \addlegendentry{$p=2$}
            \plotlinestylethree table[x index=1, y index=2, col sep=space, skip first n=1]
                {results/convergence/sbm_p3.txt};
            \addlegendentry{$p=3$}

            \plotlinestyleone[dashed]
            table [x index=1, y index=2, col sep=space, skip first n=1]
                {results/convergence/cutFEM_p1.txt};
            \plotlinestyletwo [dashed]
            table[x index=1, y index=2, col sep=space, skip first n=1]
                {results/convergence/cutFEM_p2.txt};
            \plotlinestylethree [dashed]
            table[x index=1, y index=2, col sep=space, skip first n=1]
                {results/convergence/cutFEM_p3.txt};

            \addplot[dotted, black, thick] coordinates {(0.25, 1e-2)  (0.01, 1.6e-5)} node[pos=0.95,
                    anchor=south west, yshift=0.1em]
                {$h^2$};
            \addplot[dotted, black, thick  ] coordinates {(0.25, 4* 5e-4) (0.01, 4* 3.2e-8)}
            node[pos=0.95, anchor=south west, yshift=0.2em]
                {$h^3$};

            \addplot[dotted, black, thick] coordinates {(0.25,4* 4e-5)	(0.01, 4* 1.024e-10)}
            node[pos=0.95, anchor=south west, yshift=0.3em]
                {$h^4$};
        \end{axis}
    \end{tikzpicture}
    \caption{Convergence of the $L^2$ error for the Poisson problem on a unit circle with a manufactured solution. Both
        SBM (solid lines) and CutFEM (dashed lines) demonstrate optimal convergence rates of $O(h^{p+1})$ for
        polynomial degrees
        $p=1, 2, 3$. The optimal rates are illustrated with black dotted lines. \label{fig:convergence_comparison}}

\end{figure}
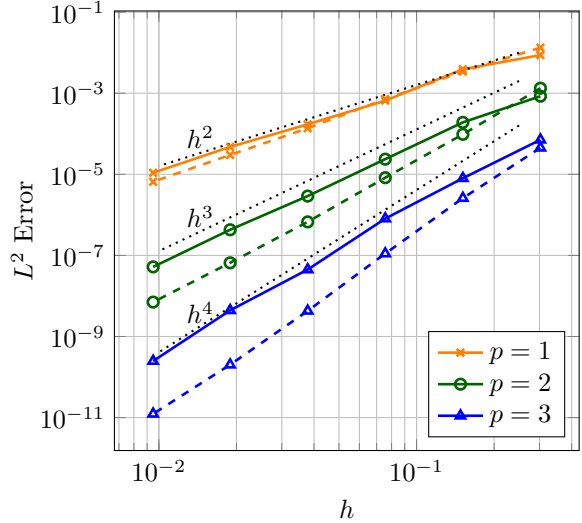

Figure~\ref{fig:convergence_comparison} displays the convergence plots for both SBM and CutFEM. The $L^2$ error is
plotted
against the mesh size $h$ on a log-log scale. For reference, lines indicating the optimal convergence rate of
$O(h^{p+1})$ are also shown.

The results confirm that both the SBM and CutFEM implementations achieve the theoretically expected
optimal convergence rates of $O(h^{p+1})$ for polynomial degree $p$. For this particular test case, CutFEM achieves
slightly smaller errors for the same mesh size $h$, which is consistent with the fact that CutFEM integrates over
the exact domain $\Omega$ while SBM uses the surrogate domain $\tilde{\Omega}$ with extrapolated boundary
conditions.
It is worth noting that the accuracy of SBM can be further improved by including some intersected cells into
the surrogate domain, as demonstrated in~\cite{yang2024optimal}.
The comparable accuracy, combined with the computational trade-offs discussed in the main body of this
paper, provides a more complete picture for choosing between the two methods.

\end{document}